\documentclass[10pt,reqno]{amsart}
\usepackage{amssymb,amsmath,amsthm}

\newcommand{\n}{\noindent}

\newcommand{\R}{\mathbb R}

\newcommand{\ZZ}{{\bar Z}}

\newcommand{\C}{\mathbb C}

\DeclareMathOperator{\Imm}{Im}
\DeclareMathOperator{\supp}{supp}

\DeclareMathOperator{\Rre}{Re}
\DeclareMathOperator{\Dom}{Dom}

\newcommand{\p}{\partial}

\newcommand{\les}{\lesssim}

\newcommand{\z}{\bar z}
\newcommand{\w}{\bar w}
\newcommand{\dbar}{\bar\partial}
\newcommand{\vp}{\varphi}

\newcommand{\atopp}[2]{\genfrac{}{}{0pt}{2}{#1}{#2}}
\newcommand{\tK}{\tilde k}
\newcommand{\tT}{\tilde T}
\newcommand{\A}[2]{a_{#1}^{#2}}

\newcommand{\nn}{\nonumber}

\newcommand{\Ts}{T_\tau}
\newcommand{\Ks}{K_\tau}
\newcommand{\LL}{\bar L}
\newcommand{\Zs}{Z_p}
\newcommand{\Zbs}{\bar{Z}_p}
\newcommand{\Zstp}{Z_{\tau p}}
\newcommand{\Zbstp}{\bar{Z}_{\tau p}}
\newcommand{\Zstpw}{Z_{\tau p,w}}

\newcommand{\Zstpz}{Z_{\tau p,z}}
\newcommand{\Zbstpz}{\bar{Z}_{\tau p,z}}
\newcommand{\Ms}{M_{\tau p}}

\newcommand{\Wstp}{W_{\tau p}}
\newcommand{\Wbstp}{\overline{W}_{\tau p}}
\newcommand{\Wstpw}{W_{\tau p,w}}
\newcommand{\Wbstpw}{\overline{W}_{\tau p,w}}

\newcommand{\Kse}{K_{\tau,\epsilon}}
\newcommand{\Ke}{k_{\ep}}
\newcommand{\vps}{{\hat\varphi}}
\newcommand{\Ystp}{Y_{\tau p}}
\newcommand{\Xstp}{X_{\tau p}}
\newcommand{\epl}{e^{i \tau T(w,z)}}
\newcommand{\emi}{e^{-i \tau T(w,z)}}
\newcommand{\T}{\epl \frac{\p}{\p\tau} \emi}

\newcommand{\ep}{\epsilon}

\newcommand{\G}{\tilde G}

\newcommand{\Gtp}{\G_{\tau p}}
\newcommand{\Boxtp}{\Box_{\tau p}}

\newcommand{\XX}{\mathcal{X}}
\newcommand{\YY}{\mathcal{Y}}
\newcommand{\MM}{\mathcal{M}}
\newcommand{\cL}{\mathcal{L}}
\newcommand{\cLb}{\bar{\mathcal{L}}}

\newcommand{\diam}[1]{\text{diam}(#1)}

\newcommand{\dni}{d_{NI}}

\newtheorem{Theorem}{Theorem}[section]
\newtheorem{Proposition}[Theorem]{Proposition}
\newtheorem{Lemma}[Theorem]{Lemma}
\newtheorem{Corollary}[Theorem]{Corollary}
\newtheorem{Remark}[Theorem]{Remark}

\newtheorem{Definition}[Theorem]{Definition}

\newcommand{\toap}{\textit{to appear,} J. Geom. Anal.}
\numberwithin{equation}{section}

\begin{document}


\author{Andrew S. Raich}

\title[OPF Operators]
{{\large One-Parameter Families of Operators in {$\C$}}
\footnote{\toap, \textbf{16}(2):353-374, 2006.} }




\maketitle
\begin{abstract}
We introduce classes of one-parameter families (OPF) of operators on $C^\infty_c(\C)$ which 
characterize the behavior of operators associated to the $\dbar$-problem in the weighted space
$L^2(\C,e^{-2p})$ where $p$ is a subharmonic, nonharmonic polynomial. 
We prove that an order 0 OPF operator extends to a bounded
operator from $L^q(\C)$ to itself,
$1<q<\infty$, with  a bound that depends on $q$ and the degree of $p$ but not on the
parameter $\tau$ or
the coefficients of $p$. Last, we show that there is a one-to-one correspondence
given by the partial Fourier transform in $\tau$ between OPF operators of order
$m\leq 2$ and nonisotropic smoothing (NIS) operators of order $m\leq 2$ on
polynomial models in $\C^2$. 
\end{abstract}

%
%
\section{Introduction.}\label{sec:intro}

The goal of this paper is to introduce classes of one-parameter families (OPF) of operators
on $\C$
which characterize the behavior of kernels associated to the weighted $\dbar$-problem in
$\C$. The need for OPF operators stems from problems associated to the inhomogeneous
$\dbar_b$-equation on polynomial models in $\C^2$  and the
$\dbar$-problem in weighted $L^2$ spaces in  $\C$. 
A polynomial model  $M$ is the boundary of an unbounded weakly pseudoconvex domain of 
finite type of the form 
$\{(z_1,z_2)\in\C^2 : \Imm z_2 > p(z_1)\}$
where $p$ is a subharmonic, nonharmonic polynomial. 
$M\cong \C\times \R$ and $\dbar_b$, defined
on $M$, can be identified with the vector field $\bar L = \frac{\p}{\p\z} - i \frac{\p p}{\p\z}\frac{\p}{\p t}$.
Under the partial Fourier transform in $t$, the vector field $\bar L$ becomes 
\begin{equation}\label{eq:zbstp} 
 \Zbstp =  \frac{\p}{\p\z} + \tau\frac{\p p}{\p\z},
\end{equation}
which we regard as 
a one-parameter family of differential operators acting on functions defined on  $\C$.
OPF operators will be defined so that
$\Zbstp$ and $\Zstp = -\Zbstp^* = \frac{\p}{\p z} - \tau\frac{\p p}{\p z}$ are the natural differential operators
under whose action OPF operators behave well.  

When $\tau=1$, the differential operator 
$\Zbs =  \frac{\p}{\p\z} + \frac{\p p}{\p\z}$ has been well studied \cite{Christ91,Ber96,Rai05h,Rai06}.
Christ \cite{Christ91} and the author \cite{Rai05h, Rai06} 
expressly cite the study of $\dbar_b$ on polynomial models as motivation to study 
the $\dbar$-problem on weighted $L^2$ in $\C$. 
In Section \ref{subsec:connection}, we review the equivalence of  the 
$\dbar$-problem in $L^2(\C, e^{-2p})$ with
the $\Zbs$-problem, $\Zbs u = f$, in $L^2(\C)$.
When $p$ is a subharmonic function satisfying  mild hypotheses on $\triangle p$,
Christ \cite{Christ91} solves the equation $\Zbs u =f$ on $L^2(\C)$  via
the complex Green operator $G_p$ 
for $\Box_p = -\Zbs \Zs$ where $\Zs = -\Zbs = \frac{\p}{\p z} - \frac{\p p}{\p z}$.
Both $G_p$ and the relative fundamental solution  $Z_p G_p$ are given as fractional integral
operators. Also, Christ shows that if $Y^\alpha$ is a product of length 2 of operators 
of the form $Y = \Zbs$ of $\Zs$,
then $Y^\alpha G_p$ is bounded on $L^q(\C)$, $1<q<\infty$. When $\tau=1$,
$G_p$ serves as a 
model for an order 2 OPF operator, while $Y^\alpha G_p$ serves as model for an
order 0 OPF operator. Christ and the author \cite{Rai05h} find pointwise estimates of the integral
kernel of  $\Gtp$ and its derivatives (Christ in the case $\tau=1$ and the author for $ \tau>0$),
and the author \cite{Rai06} finds cancellation conditions for $G_{\tau p}$ and its derivatives when $\tau>0$.
Similarly to the ordinary Laplace operator, $\Box_p$ is a 
second order, nonnegative elliptic operator, and
there is a strong analogy between $G_p$ and the Newtonian potential $N$
on $\C$. Both invert ``Laplace" operators, and if $D^2$ is a second order derivative, 
$D^2 N$ is a Calder\`on-Zygmund operator
and bounded on $L^q$, $1<q<\infty$. 
In Theorem \ref{thm:Lp bound}, we will see that order
0 OPF operator is bounded in $L^q$, $1<q<\infty$.

\subsection{Connection of $\Zbstp$ with $\dbar u =f$ on weighted $L^2$.}\label{subsec:connection}

H\"ormander's work \cite{Hor65} on solving the inhomogeneous Cauchy-Riemann equations
on pseudoconvex domains in $\C^n$. H\"ormander's methods, now classical in the subject \cite{Hor90},
rely on proving that if $\diam\Omega\leq1$, there is a solution to $\dbar u=f$ satisfying
in $L^2(\Omega,e^{-2p})$ satisfying the estimate
$\int_{\Omega} |u|^2 e^{-2p}\, dz \leq  \int_{\Omega} |f|^2 e^{-2p}\, dz$. Using the techniques 
of H\"ormander, Forn{\ae}ss and Sibony \cite{FoSi91} generalize the $L^2$ estimate to an 
$L^q$ estimate, $1<q\leq 2$, and
prove that $\dbar u=f$ has a solution satisfying:
$\left(\int_{\Omega} |u|^q e^{-2p}\, dz\right)^{\frac 1q} 
\leq \frac C{p-1}\left(\int_{\Omega} |f|^q e^{-2p}\, dz\right)^{\frac 1q}$. 
They also show that the estimate fails if $q>2$. Berndtsson \cite{Ber92} builds on the 
work of Forn{\ae}ss and Sibony and shows an $L^q$-$L^1$ result. He shows that if
$\text{diam}\,\Omega <1$ and $1\leq q < 2$, then $\dbar u=f$ has a solution so that
$\|u e^{-p}\|_{L^q(\Omega)} \leq C_q \|f e^{-p} \|_{L^1(\Omega)}$.
Berndtsson also proves a weighted $L^\infty$-$L^q$ estimate when $q>2$.

In \cite{Christ91}, Christ recognizes that it is possible to study the $\dbar$-problem in 
$L^2(\C,e^{-2p})$ by working with a related operator in the unweighted space $L^2(\C)$.
If $\dbar\tilde u =\tilde  f$ 
and both $\tilde u = e^{p} u$ and $\tilde  f= e^pf$ are in $L^2(\C,e^{-2 p})$, then
$\frac{\p\tilde u}{\p\z} = \tilde f \Longleftrightarrow e^{-p} \frac{\p}{\p\z} e^p u =f$.
However, $e^{-p} \frac{\p}{\p\z} e^p u = \ZZ_p  u$.
Consequently, the $\dbar$-problem on $L^2(\C,e^{-2p})$ is equivalent to
the $\ZZ_p $-problem, $\ZZ_p  u=f$, on $L^2(\C)$. Berndtsson \cite{Ber96} solves
$\Zbs u=f$ on smoothly bounded domains in $\C$ and views $\Box_p$ from the viewpoint
of mathematical physics.
He writes $\Box_p$ as a magnetic Schr\"odinger operator with an electric potential and his estimates
follow from Kato's inequality, a result from mathematical physics. The author \cite{Rai05h} solves
the heat equation associated to $\Box_p$ and uses techniques both from mathematical physics
and the solution of the $\Box_b$-heat equation on polynomial models in $\C^2$ \cite{NaSt00}.

\subsection{The relationship between NIS and OPF operators.}\label{subsec:nisopf}

For computations involving $\dbar_b$ on both polynomial models in $\C^2$ and the 
boundaries of other
weakly pseudoconvex domains of finite type in $\C^n$, nonisotropic smoothing operators
(NIS) operators have played a critical role in the analysis of the relative fundamental
solutions of $\Box_b$ and related operators. 
Nagel et al.\ \cite{NaRoStWa89} introduce NIS operators
while analyzing of the Szeg\"o kernel on
weakly pseudoconvex domains of finite type in $\C^2$. Nagel and Stein use properties of
NIS operators in their analysis of the heat kernel on
polynomial models in $\C^2$ \cite{NaSt00} and both the
relative fundamental solution of $\Box_b$ and the Szeg\"o kernel
on product domains and decoupled domains in $\C^n$ \cite{NaSt04,NaSt03}. 
A motivation for developing NIS operators is that the class of NIS operators
have invariances that individual operators do not. NIS operators are invariant
under translations and dilations, derivatives of NIS operators are again an NIS
operators, and order 0 NIS operators have desirable mapping properties, namely 
$L^p$-boundedness \cite{NaSt04}.

In \cite{NaSt00}, Nagel and Stein solve the $\Box_b$-heat equation
$\frac{\p u}{\p s} + \Box_b u=0$ with initial condition $u(0,\alpha) = f(\alpha)$ where
$s\in (0,\infty)$ and $\alpha\in\C\times\R$. They write their solution using the heat semigroup
$e^{-s\Box_b}$ and in turn express $e^{-s\Box_b}[f]$ as integration against a kernel called
the heat kernel. NIS operators are one of the
workhorses of their arguments because as a class of operators,
NIS operators (1) commute with vector fields
$\bar L$ and $\bar L^*$,
(2) remain invariant under translations and scaling, and 
(3) change products of arbitrary compositions of $\bar L$ and $\bar L^*$ to a composition
of a power of $\Box_b$ with a well-controlled NIS operator. The analogy of NIS
operators with Calder\`on-Zygmund operators is
strong. For example,  (3) is analogous to writing an arbitrary derivative as the composition
of $\triangle^k$
for some $k$ with a Riesz transform.

A goal for OPF operators is to play the analogous role for 
objects associated to the operators $\Zbstp$ and $\Zstp$
as NIS operators do to 
objects related to $\dbar_b$ and $\dbar_b^*$ defined on the boundaries of
weakly pseudoconvex domains in $\C^2$.
In \cite{Rai05h, Rai06}, the author solves the
$\Boxtp$-heat equation  for $\tau\in\R$, i.e. he solves the equation
$\frac{\p u}{\p s}+\Boxtp u=0$ with initial condition $u(0,z) =f(z)$. The solution is written
as integration against a kernel, called the heat kernel which is shown to be smooth off
of the diagonal $\{(s,z,w):s=0 \text{ and }z=w\}$. Also, the author finds pointwise
decay estimates for the heat kernel and its derivatives.
OPF operators play a fundamental role in these articles. They are an essential tool
in the regularity arguments and the derivative estimates.
Also, the ability to scale an OPF operator
and stay withn the class of OPF operators is crucial in the time decay estimate of the heat
kernel $e^{-s\Boxtp}$.

\section{Main Results.}\label{sec:results}


\begin{Theorem}\label{thm:Lp bound}
If $\Ts$ is an OPF operator of order 0, then $\Ts$, $\Ts^*$ are bounded operators from
$L^q(\C)$ to $L^q(\C)$, $1<q<\infty$, with a constant independent of $\tau$ but
depending on $q$.
\end{Theorem}
Also, the classes of OPF operators fulfill the 
promise of being an analog to NIS operators. We can use results about OPF operators
to study NIS operators and vice versa. We have the theorem:
%
%
\begin{Theorem}\label{thm:correspondence}
Given a subharmonic, nonharmonic polynomial $p:\C\to\R$, 
there is a one-to-one correspondence between OPF operators
of order $m\leq 2$ with respect to $p$
and NIS operators of order $m\leq 2$ on the polynomial model
$M^p = \{(z_1,z_2)\in\C^2 : \Imm z_2 = p(z_1)\}$. The correspondence is given by a partial
Fourier transform in $\Rre z_2$.
\end{Theorem}

%
%
%
\section{Notation and Definitions.}\label{sec:notation}

\subsection{Notation For Operators on $\C$.}\label{subsec:notationdown}
For the remainder of the article,  let $p$ be a  subharmonic, nonharmonic polynomial.
It will be important for us to expand $p$ around an arbitrary point $z\in\C$, and 
we set:
\begin{equation}\label{eq:Ajk}
\A{jk}z = \frac{1}{j! k!}\frac{\p^{j+k}p}{\p z^j\p\z^k}(z).
\end{equation} 
We need the following two ``size" functions to write down the size and cancellation conditions 
for both OPF operators and NIS operators. Let
\begin{align}\label{eq:Lam}
\Lambda(z,\delta) &= \sum_{j,k\geq 1}  \left|\A{jk}z \right|\delta^{j+k}
\intertext{and}
\label{eq:mu}
\mu(z,\delta) &= \inf_{j,k\geq 1} \frac{|\delta|^{1/j+k}}{|\A{jk}z|^{1/j+k}}.
\end{align}
It follows 
$\mu(z,\delta)$ is an approximate inverse to $\Lambda(z,\delta)$. This means that if $\delta>0$,
\begin{equation}\label{eqn:mu,Lambda approx inverses}
\mu\big(z,  \Lambda(z,\delta)\big) \sim \delta
\quad\text{ and }\quad
\Lambda\big(z,\mu(z,\delta)\big) \sim \delta.
\end{equation}
We use the notation $a\les b$ if $a\leq C b$ where $C$ is a constant that may depend on the
dimension 2 and the degree of $p$. We say that $a\sim b$ if $a\les b$ and $b\les a$.

$\Lambda(z,\delta)$ and $\mu(z,\delta)$ are geometric objects from the Carnot-Carath\'eodory
geometry  developed by Nagel et al.\  
\cite{NaStWa85,Na86}. The functions also arise
in the analysis of  magnetic Schr\"odinger  operators with electric potentials \cite{Shen96,
Shen99,Kurata00, Rai05h,Rai06}.

Denote the ``twist" at $w$, centered as $z$ by 
\begin{align}
T(w,z) &= - 2\Imm \left( \sum_{j \geq 1} \frac{1}{j!} \frac{\p^j p}{\p z^j}(z) (w-z)^j\right) \nn\\
&= i \left(  \sum_{j \geq 1} \frac{1}{j!} \frac{\p^j p}{\p z^j}(z) (w-z)^j
  -  \sum_{j \geq 1} \frac{1}{j!} \frac{\p^j p}{\p \z^j}(z) \overline{(w-z)}^j\right).
\label{eq:twist}
\end{align}

Also associated to a polynomial $p$ and the parameter $\tau\in\R$ are the weighted differential
operators 
\[
\Zbstpz = \frac{\p}{\p \z} + \tau \frac{\p p}{\p \z}= e^{-\tau p}\frac{\p p}{\p \z}e^{\tau p}
\hspace{1cm}
\Zstpz = \frac{\p}{\p z} - \tau \frac{\p p}{\p z}  = e^{\tau p}\frac{\p p}{\p z}e^{-\tau p}.
\]
We need to establish notation for adjoints. 
If $T$ is an operator (either 
bounded or closed and densely defined) on a Hilbert space with inner product
$\big(\,\cdot\, ,\cdot\,\big)$, let $T^*$ be the Hilbert space adjoint of $T$. This means that if
$f\in\Dom T$ and $g\in \Dom{T^*}$, then $\big(Tf,g\big) = \big(f,T^*g\big)$. The Hilbert spaces that
arise in this paper are $L^2(\C)$ and $L^2(\C\times\R)$.
Since the $L^2$-adjoints of $\Zbstp$ and $\Zstp$ are different than their adjoints in the sense of
distributions, for clarity we let $\Wbstp$ and $\Wstp$ be the negative of 
the distributional adjoints of $\Zbstp$ and
$\Zstp$, respectively. Thus,
\[
\Wbstpw = \frac{\p}{\p \w} -  \tau \frac{\p p}{\p \w}= e^{\tau p}\frac{\p p}{\p \w}e^{-\tau p}
\hspace{1cm}
\Wstpw  = \frac{\p}{\p w} +  \tau \frac{\p p}{\p w}= e^{-\tau p}\frac{\p p}{\p w}e^{\tau p}.
\]
We think of $\tau$ as fixed and the operators $\Zbstpz$, $\Zstpz$, $\Wbstpw$, and $\Wstpw$ as acting
on functions defined on $\C$.
Also, we will omit the variables $z$ and $w$ from 
subscripts when the application is unambiguous. Observe that
$\overline{(\Zstp)} =  \Wbstp$ and $\overline{(\Zbstp)} = \Wstp$.
Finally, let
\[
 \Ms = \T.
\]

\subsection{Definition of OPF Operators.}\label{subsec:OPF operator}
\renewcommand{\labelenumi}{(\alph{enumi})}
Let $p$ be a subharmonic, nonharmonic polynomial.
We say that  $\Ts$ is a \emph{one-parameter family (OPF)
of operators}  of order $m$ with respect to the polynomial $p$ if the 
following conditions hold:
\begin{enumerate}  

\item There is a  function 
$\Ks\in C^\infty\Big(\big((\C\times\C) \setminus\{z=w\}\big) \times (\R\setminus\{0\})\Big)$ 
so that for fixed $\tau$,  
$\Ks$ is a distributional kernel, i.e. if $\vp,\psi \in C^\infty_c(\C)$ and 
$\supp\vp \cap \supp \psi = \emptyset$, then
$\Ts[\vp]\in (C^\infty_c)'(\C)$ and  
\[
 \langle \Ts[\vp](\cdot),\psi\rangle_{\C} = \iint_{\C\times\C} \Ks(z,w)\vp(w) \psi(z)\, dw dz.
\]

\item \label{it:approx}
There exists a family of functions $\Kse(z,w)\in C^{\infty}(\C\times\C\times\R)$ so
that if $\vp\in C^\infty_c(\C\times\R)$,
\[
 \Kse[\vp]_{\C\times\R}(z,\tau) = \int_{\C\times\R}\vp(w,\tau)\Kse(z,w)\,dw d\tau
\]
and $\lim_{\epsilon\to0} \Kse[\vp]_{\C\times\R}(z) = \Ks[\vp]_{\C\times\R}(z)$ in 
$(C^\infty_c)'(\C\times\R)$.  \\
All of the additional conditions are assumed to apply to the kernels 
$\Kse(z,w)$ uniformly in $\ep$.

\item 
Size Estimates. If $\Ystp^{J}$ is a product of $|J|$ operators of the form
$\Ystp^j =\Zstpz$, $\Zbstpz$,  $\Wstpw$, $\Wbstpw$, or
$\Ms$ where $|J| = \ell +n$ and
$n = \#\{j : \Ystp^j = \Ms\}$, for any $k\geq 0$ there exists a constant $C_{\ell,n,k}$ so that
\begin{equation}\label{it:size}
\left| \Ystp^{J}  \Kse(z,w)\right| \leq C_{\ell,n,k} 
\frac{|z-w|^{m-2-\ell}}{|\tau|^{n+k}\Lambda(z,|w-z|)^k} \quad \text{if}\quad
\begin{cases} &m<2 \\ &m=2,\ k\geq 1 \\ &m=2, |w-z| > \mu(z,\tfrac 1\tau) \end{cases} 
\end{equation}

Also, if $m=2$ and $|w-z| \leq \mu(z,\tfrac 1\tau)$, then
\begin{equation}\label{it:size2}
\left| \Ms^n \Kse(z,w) \right| \leq C_n 
\begin{cases} \log\left( \frac{2\mu(z,\tfrac 1\tau)}{|w-z|}\right)  &n=0 \\
|\tau|^{-n} &n\geq 1 \end{cases}
\end{equation}

\item 
Cancellation in $w$.  If $\Ystp^{J}$ is a product of $|J|$ operators of the form
$\Ystp^j = \Zstpz$, $\Zbstpz$, $\Wbstpw$, $\Wstpw$, or $\Ms$ where $|J| = \ell +n$ and
$n = \#\{j : \Ystp^j = \Ms\}$, for any $k\geq 0$ there exists a constant $C_{\ell,n,k}$ and $N_\ell$ 
so that for $\vp\in C^\infty_c (D(z_0,\delta))$,
\begin{align}\label{it:cancel w}
& \sup_{z\in\C} \left| \int_{\C}  \Ystp^{J}\Kse(z,w)\vp(w)\,dw\right|\nn \\
& \leq \frac{C_{\ell,n,k}}{|\tau|^n} 
 \begin{cases} 
   {\displaystyle\delta^2
  \Big(\log\big(\tfrac{2\mu(z,\frac 1\tau)}{\delta}\big) \|\vp\|_{L^\infty(\C)}
  + \sum_{1 \leq |I|\leq N_0} \delta^{|I|}\|\Xstp^I\vp(w)\|_{L^\infty(\C)}\Big)} 
  &\atopp{\displaystyle \delta < \mu(z,\tfrac 1\tau) 
\text{ and}}{\displaystyle m=2,\ell=0} \\
 {\displaystyle \frac{\delta^{m-\ell}}{|\tau|^{k}\Lambda(z,\delta)^k} \sum_{|I| \leq N_\ell}
 \delta^{|I|} \left\| \Xstp^I \vp\right\|_{L^\infty(\C)}} & \text{otherwise}
\end{cases}
\end{align}

 where $\Xstp^I$ is composed solely of  $\Zstp$ and $\Zbstp$.

\item Cancellation in $\tau$. If $\Xstp^{J}$ is a product of $|J|$ operators of the form
$\Xstp^j = \Zstpz,\ \Zbstpz$ or $\Wstpw,\ \Wbstpw$ 
and $|J| = n$,
there exists a constant $C_{n}$ so that
\begin{equation} \label{it:cancel tau}
 \int_\R\Xstp^{J}\left( e^{ i \tau t}\Kse(z,w)\right) \, d\tau \leq 
 C_n \frac{\mu(z,t+T(w,z))^{m-n}}{\mu(z,t+T(w,z))^2 |t+T(w,z)|}.
\end{equation}

\item \label{it:adjoint} Adjoint.
Properties (a)-(e) 
also hold for the adjoint operator $\Ts^*$ whose distribution kernel is given
by $\overline{\Kse(w,z)}$
\end{enumerate}

Note that for  the $\tau$-cancellation condition \eqref{it:cancel tau}, we do not need to
consider the case $\Xstp^j = \Ms$ since 
$\int_\R \frac{\p}{\p\tau}\big( e^{ i \tau (t+T(w,z))}\Kse(z,w)\big) \, d\tau =0$.

In the size condition (c) and cancellation condition (d), the 
$\tau^k \Lambda(z,|z-w|)^k$ and $\tau^k \Lambda(z,\delta)^k$ terms indicate rapid
decay. If OPF operators are to be partial Fourier transforms of NIS operators
on polynomial models, rapid decay should not be surprising;
it is consequence of being able to integrate
parts from the Fourier transform formula. This will be seen explicitly in Lemma \ref{lem:size up to down}.
Ignoring the rapid decay terms, the size and cancellation conditions of OPF operators 
are familiar. An order 2 OPF operator should ``invert" two derivatives, like the 
Newtonian potential. In $\R^2$, the Newtonian potential has a logarithmic
blowup on the diagonal, just like an order 2 OPF operator. For an order 0 OPF operator,
the blowup on the diagonal 
is the same as a Calder\`on-Zygmund kernel, and the decay of
$\Ks(0,z)$ is $|z|^{-2}$, the same as a Calder\`on-Zygmund kernel.
For the cancellation conditions, if
$\vp$ is ``normalized" appropriately, the cancellation condition \eqref {it:cancel w} simplifies to
\[
\|\Ystp^J \Ts[\vp]\|_{L^\infty(\C)} \les \delta^j.
\]
This is reminiscent of cancellation of a Calder\`on-Zygmund operator or an NIS operator.

\subsection{Notation for Carnot-Carath\'eodory geometry and Vector Fields on $\C\times\R$.}\label{subsec:CC}
In order to write down the definition of an NIS operator
on a polynomial model  in $\C^2$, we need to establish notation
for the Carnot-Carath\'eodory metric $\rho$ and corresponding balls  $B_{NI}\big((z,t),\delta\big)$. 
If $M^p$ is a polynomial model  in $\C^2$ given by
$M^p = \{(z_1,z_2)\in\C^2 : \Imm z_2 = p(z_1)\}$, then $M^p \cong \C\times\R$. Under
the isomorphism, a representation of the  Carnot-Carath\'eodory  metric 
is the nonisotropic pseudodistance $\rho\big((z,t),(w,s)\big) = |z-w| + \mu\big(z,t-s+T(w,z)\big)$
where $(z,t),(w,s)\in\C\times\R$.
Since $\rho\big((z,t),(w,s)\big)$ is a function of $z$, $w$, and $t-s$, we define a new function
\begin{equation}\label{eq:dni}
\dni(z,w,t) = |z-w| + \mu\big(z,t+T(w,z)\big).
\end{equation}
We will see that $\dni(z,w,t)$ is essentially symmetric in $(z,w)$.
The nonisotropic ball 
\[
B_{NI}\big((z,t),\delta\big) = \{ (w,s)  :  \dni(z,w,t-s) <\delta \}.
\]
We also define a volume function
\[
V_{NI}\big((z,t),(w,s)\big) = \big|B_{NI}\big((z,t),\dni(z,w,t-s)\big)\big|
\sim \dni(z,w,t-s)^2 \Lambda\big(z,\dni(z,w,t-s)\big).
\]
That the volume function is comparable to $\dni(z,w,t-s)^2 \Lambda\big(z,\dni(z,w,t-s)\big)$ follows from
\eqref{eqn:mu,Lambda approx inverses}.

If $\tau$ is the transform variable of $t$, 
observe that under the partial Fourier transform in $t$, 
$\Zbstp$ and $\Zstp$ map to the vector fields
\begin{align*}
&\LL_z = \frac{\p}{\p \z} - i \frac{\p p}{\p \z} \frac{\p}{\p t}&  
&L_z = \frac{\p}{\p z} + i \frac{\p p}{\p z} \frac{\p}{\p t} 
\intertext{while  $\Wbstp$ and $\Wstp$ map to the vector fields}
& \cLb_w = \frac{\p}{\p \w} + i\frac{\p p}{\p\w}\frac{\p}{\p t}&
& \cL_w = \frac{\p}{\p w} - i\frac{\p p}{\p w}\frac{\p}{\p t}.
\end{align*}
As we know from Section \ref{sec:intro},
$\dbar_b$ (defined on $M$) becomes the operator $\LL_z$ on $\C\times\R$. It follows 
that $-L_z$ is the Hilbert space adjoint to $\LL_z$ in $L^2(\C\times\R)$. The translation invariance
in $t$ causes many operators of interest to have a convolution structure in $t$. A consequence
is that if we have a function $\tilde f\big((z,t),(w,s)\big) = f(z,w,t-s)$, we may study
$f(z,w,t)$. By the chain rule, $\cLb_w$ and $\cL_w$ are  the versions of $\LL_z$
and $L_z$ in the $w$-variable. 
Finally, let
\[
\MM = - i\big(t+T(w,z)\big).
\]

\subsection{NIS operators on polynomial models in $\C^2$.} \label{subsec:NIS}
There are different notions of NIS operators (e.g. \cite{NaRoStWa89,NaSt00}). We use the
definition from \cite{NaRoStWa89}. 

\begin{Definition}[\hspace{-4pt}. Nonisotropic Smoothing Operator of order $m$]\label{def:NIS}
\renewcommand{\labelenumi}{(\alph{enumi})}

Let 
\[ 
T[f](z,t) = \int_{\C\times\R} T\big((z,t),(w,s)\big)f(w,s)\, dw ds, 
\] 
where
$T\big((z,t),(w,s)\big)$ is a distribution which is $C^\infty$ away from the diagonal.
We shall say that $T$ is a nonisotropic smoothing operator which is smoothing of
order $m$ if
there exists a family 
\[
T_\ep[f](z,t) = \int_{\C\times\R} T_\ep\big((z,t),(w,s)\big) f(w,s)\, dwds,
\]
so that:
\begin{enumerate}
\item $T_{\ep}[f] \to T[f]$ in $C^\infty(\C\times\R)$ as $\ep\to 0$ whenever
$f\in\C^\infty_c(\C\times\R)$;
\item Each $ \displaystyle T_\ep\big((z,t),(w,s)\big) \in C^\infty\big((\C\times\R)\times(\C\times\R)\big)$;\\
The following two conditions hold uniformly in $\ep$:
\item If $\XX^I = \XX_{i_1}\XX_{i_2}\cdots \XX_{I_k}$ where
$\XX_{i_j} = L_z, L_w, \LL_z$, or $\LL_w$, then
\begin{equation}\label{eq:NIS size}
\big| \XX^I T_\ep\big((z,t),(w,s)\big) | \leq c_{|I|} \frac{\dni(z,w,t-s)^{m-|I|}}
{V\big((z,t),(w,s)\big)};
\end{equation}
\item For each $\ell\geq 0$, there exists an $N = N_\ell$ so that whenever $\vp$
is a smooth (bump) function supported in $B_{NI}\big((z,t),\delta\big)$,
\begin{equation}\label{eq:NIS cancel}
\big| \XX^IT[\vp](z,t)\big| \leq C_{\ell} \delta^{m-\ell} \sup_{w,s} \sum_{|J|\leq N_\ell}
\delta^{|J|} \big|\XX^J[\vp](w,s)\big|,
\end{equation}

 $|I|=\ell$;
\item The same estimates hold for the adjoint operator $T^*$, i.e. the operator with the kernel
$\overline{T\big((w,s),(z,t)\big)}$.
\end{enumerate}
\end{Definition}

\section{Properties of $T(w,z)$.}

To prove Theorem \ref{thm:Lp bound} and Theorem \ref{thm:correspondence}, we need to 
understand the ``twist" $T(w,z)$ and how it behaves under differentiation.
\begin{Proposition}\label{prop:T(w,z) = -T(z,w)}
\[
T(w,z) =- T(z,w).
\]
\end{Proposition}

\begin{proof}
Since $p(z) = \sum_{j,k} \frac{1}{j!k!} \frac{\p^{j+k}p}{\p z^j \p\z^k}(w) (z-w)^j
\overline{(z-w)}^k$,  we have
\begin{align*}
\frac{\p^\ell p}{\p z^{\ell}}(z) =
& \sum_{\atopp{j\geq \ell}{k\geq 0}} \frac{j!}{(j-\ell)!} \frac{1}{j!k!}.
\frac{\p^{j+k}p}{\p z^j \p\z^k}(w)(z-w)^{j-\ell} \overline{(z-w)}^k. 
\end{align*}
Since $p$ is $\R$-valued, the twist [Equation \eqref{eq:twist}]
$T(w,z) = -2\Imm \left( \sum_{\ell \geq 0} \frac{1}{\ell!} \frac{\p^\ell p}{\p z^\ell}(z) (w-z)^\ell\right)$, so
\begin{align*}
\sum_{\ell\geq 0} \frac{1}{\ell!} &\frac{\p^\ell p}{\p z^\ell}(z)(w-z)^\ell \\
&= \sum_{\ell\geq 0}  \frac{1}{\ell!} \left(\sum_{\atopp{j\geq \ell}{k\geq 0}} \frac{j!}{(j-\ell)!} \frac{1}{j!k!}.
\frac{\p^{j+k}p}{\p z^j \p\z^k}(w)(z-w)^{j-\ell} \overline{(z-w)}^k\right) (w-z)^\ell \\
&= \sum_{\atopp{j\geq 0}{k\geq 0}} \left(\sum_{\ell=0}^j \binom{j}{\ell} (-1)^{\ell}\right) \frac{1}{j!k!}
\frac{\p^{j+k}p}{\p z^j \p\z^k}(w) (z-w)^j \overline{(z-w)}^k \\
&= \sum_{k\geq 0} \frac{1}{k!} \frac{\p^k p}{\p \z^k}(w)\overline{(z-w)}^k 
= \overline{ \sum_{j\geq 0} \frac{1}{j!} \frac{\p^j p}{\p z^j}(w)(z-w)^j}.
\end{align*}
The  second to last line uses the identity $\sum_{\ell=0}^j \binom{j}{\ell}(-1)^\ell = \delta_0(j)$.
The result follows easily.
\end{proof}

\begin{Corollary} \label{cor:mu symmetry}
\[
\dni(z,w,t) \sim \dni(w,z,t).
\]
\end{Corollary}

\begin{proof} This is a well known fact (\cite{NaStWa85,Na86}), but we are in a situation where
the computations can be explicit. We sketch a proof.
If $r = |t+T(w,z)|$, it follows from from Proposition \ref{prop:T(w,z) = -T(z,w)} 
that it is enough to show that
\[
|z-w| + \mu(z,r) \sim |z-w| + \mu(w,r).
\]  
If $\mu(z,r)< |z-w|$ and $\mu(w,r)<|z-w|$, there is nothing to prove, so (without loss
of generality) assume that
$\mu(z,r)>|z-w|$. 
By expanding $p(z)$ around $w$ and $p(w)$ around $z$, it can be shown  that
$\Lambda(z,\delta) \sim\Lambda(w,\delta)$ if $\delta > |w-z|$. Thus, we see
\[
\Lambda\big(w,\mu(z,r)\big)  \sim \Lambda\big(z,\mu(z,r)\big) \sim r,
\]
and it follows that $\mu(z,r)\sim\mu(w,r)$.
\end{proof}

The next proposition contains two useful, though simple, computations.
\begin{Proposition} \label{prop:twist derivative}
\[
\frac{\p T}{\p z}(w,z) = -i \frac{\p p}{\p z}(z) - i \sum_{j\geq 1} \frac{1}{j!} \frac{\p^{j+1} p}{\p z\p\z^j}(z)
\overline{(w-z)}^j 
\]
and
\[
\frac{\p T}{\p\z}(w,z) = i \frac{\p p}{\p\z}(z) + i \sum_{j\geq 1} \frac{1}{j!} \frac{\p^{j+1} p}{\p z^j\p\z}(z)
(w-z)^j. 
\]
\end{Proposition}

\begin{proof} The proof is a short computation.
\begin{align*}
\frac{\p T}{\p z}(w,z) &= i\Bigg( \sum_{j=1}^{\deg(p)-1} \frac{1}{j!} \frac{\p^{j+1} p}{\p z^{j+1}} (z) (w-z)^j
- \sum_{j=1}^{\deg(p)} \frac{1}{(j-1)!} \frac{\p^j p}{\p z^j}(z)(w-z)^{j-1} \\
&\hspace{2.65in} - \sum_{j=1}^{\deg(p)-1} \frac{1}{j!}
\frac{\p^{j+1} p}{\p z \p\z^j} (z) \overline{(w-z)}^j \Bigg)\\
&= -i \frac{\p p}{\p z}(z) - i \sum_{j\geq 1} \frac{1}{j!} \frac{\p^{j+1} p}{\p z\p\z^j}(z)
\overline{(w-z)}^j 
\end{align*}
since the first sum cancels all but the first term of the second sum. Since $T$ is $\R$-valued,
$\frac{\p T}{\p\z}(w,z) = \overline{\frac{\p T}{\p z}(w,z)}$ which gives the result
for the second sum.
\end{proof}
A useful consequence of these calculations is
\begin{Proposition}\label{prop:t+T deriv} Let $\YY^J$ be a product of $|J|$ operators of the form  
$\YY^j = L_z, \LL_z$,  $\cL_w, \cLb_w$.
Then
\[
|\YY^J\big(t+T(w,z)\big)| \leq C_{|J|} \frac{\Lambda(z,\dni(z,w,t))}{\dni(z,w,t)^{|J|}}.
\]
\end{Proposition}
Before we prove the Proposition \ref{prop:t+T deriv}, we note that
the result would be false if we replaced $t+T(w,z)$ with $t$ or $T(w,z)$. Without both terms,
there would be uncontrolled derivatives of $p$ remaining after applying $\YY^j$.
\begin{proof} We have $L_z\big(t+T(w,z)\big) = \frac{\p  T}{\p z}(w,z) + i \frac{\p p}{\p z}(z) = 
-i \sum_{j\geq 1} \frac 1{j!} \frac{\p^{j+1} p(z)}{\p z\p\z^j}\overline{(w-z)}^j$. Similarly,
$\LL_z\big(t+T(w,z)\big) = i \sum_{j\geq 1} \frac 1{j!} \frac{\p^{j+1} p(z)}{\p z^j\p \z}(w-z)^j$. Analogous 
equalities (with $z$ and $w$ interchanged and the sign switched) hold for 
$\cL_w\big(t+T(w,z)\big)$ and $\cLb_w\big(t+T(w,z)\big)$
since
\begin{align*}
&\cL_w\big(t+T(w,z)\big) = \left( \frac{\p}{\p w}- i \frac{\p p}{\p w}\frac{\p}{\p t}\right)(t- T(z,w))
= -i \frac{\p p}{\p w}(w) - \frac{\p T}{\p w}(z,w) \\
&= - \left(i \frac{\p p}{\p w}(w) + \frac{\p T}{\p w}(z,w)\right)
= - \left(\frac{\p}{\p w} + i \frac{\p p}{\p w} \frac{\p}{\p t}\right)(t+T(z,w)) = - L_w(t+T(z,w))
\end{align*}
and $\cLb_w\big(t+T(w,z)\big)= - \LL_w(t+T(z,w))$.
But
\[
\left| \sum_{j\geq 1} \frac 1{j!} \frac{\p^{j+1} p(z)}{\p z^j\p \z}(w-z)^j \right|
\leq c_1 \frac{\Lambda(z,\dni(z,w,t))}{\dni(z,w,t)}.
\]
Higher order derivatives are easier. As we just showed, the result of applying $\YY^1$ to
$t+T(w,z)$ leaves a polynomial that is a sum of derivatives of $\triangle p$ (and hence well controlled).
There are no $t$ terms remaining, so if $j\geq 2$, applying $\YY^j$ is a matter of applying one of:
$\frac{\p}{\p\z}$, $\frac{\p}{\p z}$, $\frac{\p}{\p\w}$, $\frac{\p}{\p w}$ .
Hence, the computation is simpler, and it can be done naively, i.e. there is no need to find
any cancelling terms (which in general are absent).
\end{proof}

%
%
\section{$L^q$ boundedness of order 0 operators.}\label{sec:Lp bound}

We are now ready to begin the proof Theorem \ref{thm:Lp bound}. 
The idea is to show that $\emi \Kse$ satisfies the bounds of a Calderon-Zygmund kernel
and the operator $S_\tau$ with kernel $\emi \Kse$ is restrictly bounded. These two
facts, proven in Lemma \ref{lem:Calderon-Zygmund operators} and
Lemma \ref{lem:rest bound}, respectively, 
show $S_\tau$ satisfy the hypotheses of $T(1)$ theorem \cite{St93}. Consequently,
$S_\tau$ is a bounded operator on $L^q (\C)$. A result by
Ricci and Stein \cite{RiSt87} applies to pass from
$L^q (\C)$ boundedness of $S_\tau$ to $L^q (\C)$ boundedness of $\Ts$.

\begin{Lemma}\label{lem:Calderon-Zygmund operators} Let $\Ts$ be an OPF operator 
of order $m\leq 2$ with a family of kernel approximating functions $\Kse$. 
For $k\geq 0$, there exists $C_k$ independent of $\tau$ so that $\Kse(z,w)$ satisfies:
\begin{enumerate}
\item 
\begin{equation} \label{it:dif ineq}
\left| \nabla_{z,w}\left( \emi\Kse(z,w)\right) \right| 
\leq C_{k} \frac{|w-z|^{m-3}}{|\tau|^k\Lambda(z,|w-z|)^k}
\end{equation}
\item If $2|w-w'| \leq |w-z|$, then
\begin{equation}\label{it:xi-w}
\left|\emi \Kse(z,w) - e^{- i \tau T(w',z)}\Kse(z,w')\right| \leq C_{k} \frac{|w-w'|}
{|w-z|^{3-m}|\tau|^k\Lambda(z,|w-z|)^k} 
\end{equation}
\item If $2|z-z'| \leq |w-z|$, then 
\begin{equation}\label{it:z-z'}
\hspace*{-.1in}\left|\emi \Kse(z,w) - e^{- i \tau T(w,z')}\Kse(z',w)\right| \leq C_{k}
 \frac{|z-z'|}
{|w-z|^{3-m}|\tau|^k\Lambda(z,|w-z|)^k}
\end{equation}
\end{enumerate}
Also, the constants are uniform in $\ep$.
\end{Lemma}

\begin{proof} It is immediate from the Mean Value Theorem that 
\eqref{it:dif ineq} implies \eqref{it:xi-w} 
and \eqref{it:z-z'}.
To prove \eqref{it:dif ineq}, we use Proposition \ref{prop:twist derivative} and compute:
\begin{align*}
\epl\frac{\p}{\p z}&\left( \emi \Kse(z,w)\right) 
= - i \tau\frac{\p T}{\p z}(w,z) \Kse(z,w) + \frac{\p \Kse}{\p z}(z,w) \\
&= \frac{\p \Kse}{\p z}(z,w) - \tau\frac{\p p}{\p z}(z) \Kse(z,w) -\tau
\sum_{j\geq 1}\frac{1}{j!} \frac{\p^{j+1} p}{\p z\p\z^j}(z)\overline{(w-z)}^j \Kse(z,w).
\end{align*}
Using the size estimate \eqref{it:size},
\begin{align*}
\left| \frac{\p}{\p z} \left( \emi \Kse(z,w)\right) \right|
&\leq \Zstp\Kse(z,w) + \frac{\tau\Lambda(z,|w-z|)}{|w-z|} \Kse(z,w)\\
&\leq C_k\frac{|w-z|^{m-3}}{|\tau|^k\Lambda(z,|w-z|)^k}. 
\end{align*}
A virtually identical calculation shows
\[
\left| \frac{\p}{\p \z} \left( \emi \Kse(z,w)\right) \right|
\leq C_k\frac{|w-z|^{m-3}}{|\tau|^k\Lambda(z,|w-z|)^k} 
\]
which proves $\left|\frac{\p}{\p\z}\left( \emi \Kse(z,w)\right)\right|$ 
satisfies the bound in \eqref{it:dif ineq}.  The bounds 
for the $w$ and $\w$ derivatives, $\left|\frac{\p}{\p w}\left(\emi\Kse(z,w)\right)\right|$ 
and $\left|\frac{\p}{\p\w}\left(\emi\Kse(z,w)\right)\right|$, use 
a repetition of the calculations just performed and the identity 
$\emi = e^{ i \tau T(z,w)}$ (which follows from Proposition \ref{prop:T(w,z) = -T(z,w)}).
\end{proof}

We now restrict ourselves to the case $m=0$. Given an family $\Ts$ of order 0, define a
related family of operators $S_\tau$ so 
that if $\Ks(z,w)$ is the kernel of $\Ts$, the kernel of
$S_\tau$ is given by $\emi \Ks(z,w)$.  We have the following:
\begin{Lemma} \label{lem:rest bound}
$S_\tau$ and $S_\tau^*$ are restrictly bounded uniformly in $\tau$, 
i.e. if $\vp\in C^\infty_c(D(0,1))$, $\|\vp\|_{C_{N_0}}\leq 1$
{\normalfont [}where $N_0$ is the constant from the 
cancellation condition \eqref{it:cancel w}{\normalfont ]} 
and $\vp^{R,z_0}(z) = \vp(\tfrac{z-z_0}R)$, then
\begin{equation}\label{eqn:restr bound}
 \|S_\tau (\vp^{R,z_0})\|_{L^2(\C)} \leq A R, 
\qquad \|(S_\tau)^* (\vp^{R,z_0})\|_{L^2(\C)} \leq A R
\end{equation}
with the constant $A$ independent of $\tau$.
\end{Lemma}

\begin{proof} From the adjoint condition (f), it follows that we only have the prove
the restricted boundedness of $S_\tau$. 
\begin{align*}
\| S_{\tau,\ep} (\vp^{R,z_0})\|_{L^2} & = \left( \int_{\C} \left| \int_{\C}\emi \Kse(z,w)
\vp(\tfrac{w-z_0}R)\,dw
\right|^2 dz \right)^{\frac 12} \\
&\leq  \left( \int_{|z-z_0|< 2R} \left| \int_{\C} \Kse(z,w) \left(\emi\vp(\tfrac{w-z_0}R)\right)\,dw
\right|^2 dz \right)^{\frac 12}  \\
&+  \left( \int_{|z-z_0|\geq 2R} \left| \int_{\C}\emi \Kse(z,w) \vp(\tfrac{w-z_0}R)\,dw
\right|^2 dz \right)^{\frac 12} \\
& = I + II.
\end{align*}
We estimate $I$ first. By the cancellation condition \eqref{it:cancel w}
\begin{multline*}
\left| \int_{\C} \Kse(z,w) \left(\emi\vp(\tfrac{w-z_0}R)\right)\,dw \right| \\
\leq 
C_{N_0} \frac{1}{\max\{1, |\tau|^{N_0}\Lambda(z,R)^{N_0}\}} \sup_{w\in\C}  \sum_{|I|\leq N_0}
R^{|I|} \left| \Ystp^I \left(e^{- i \tau T(w,z)}\vp(\tfrac{w-z_0}R)\right)\right|.
\end{multline*}
We claim $R^{|I|} \left| \Ystp^I \left(e^{ i \tau T(z,w)}\vp(\tfrac{w-z_0}R)\right)\right| \leq C_{|I|}
\max\{ 1, |\tau|^{|I|}\Lambda(z,R)^{|I|}\}$. To see this, we first do the case
$\Ystp^I = \Zstpw$.
It follows from Proposition \ref{prop:T(w,z) = -T(z,w)} and Proposition \ref{prop:twist derivative} that
\begin{multline*}
\Zstpw  \left(e^{ i \tau T(z,w)}\vp(\tfrac{w-z_0}R)\right)
= \frac{e^{ i \tau T(z,w)}}R \frac{\p\vp}{\p w}(\tfrac{w-z_0}R) \\ +  \tau e^{ i \tau T(z,w)}
\sum_{j\geq 1} \frac{1}{j!} \frac{\p^{j+1} p}{\p w\p\w^j}(w)\overline{(z-w)}^j \vp(\tfrac{w-z_0}R).
\end{multline*}
Hence, $\left| \Zstpw  \left(e^{ i \tau T(z,w)}\vp(\tfrac{w-z_0}R)\right)\right| \leq
\frac CR \left(1 + \tau\Lambda(z,R)\right)$.
Iterating this argument proves the claim. Thus, for $|z-z_0|\leq 2R$,
\[
\left| \int_{\C} \Kse(z,w) \emi\vp(\tfrac{w-z_0}R)\,dw \right| \leq C,
\]
and 
\[
I \leq C  \left( \int_{|z-z_0|< 2R}  dz \right)^{\frac 12}  \leq AR.
\]
When $|z-z_0|\geq 2R$, $|z-z_0| \sim |z-w|$ for $w \in \supp \vp(\tfrac{\cdot-z_0}R)$, so
\[
II \leq C \left( \int_{|z-z_0|\geq 2R} \frac{1}{|z-z_0|^4} 
\left(\int_{\C}  \left|\vp(\tfrac{w-z_0}R)\right|\,dw\right)^2
dz \right)^{\frac 12} \leq C R^2 \left(\int_{r>R} \frac 1{r^3}\,
dr\right)^{\frac 12} 
\leq A R.
\]
\end{proof}

The final ingredient we need to prove Theorem \ref{thm:Lp bound} is a result
by Ricci and Stein \cite{RiSt87}.
\begin{Theorem}[(Ricci-Stein)] \label{thm:ricci-stein}
In $\R^n\times\R^n$, let $K( \cdot\, ,\cdot)$ satisfy the following:
\begin{enumerate} 
 \item $K(\cdot\, ,\cdot)$ is  a $C^1$ function away from the diagonal 
 $\{(x,y)\in\R^n\times\R^n : x=y\}$, 
 \item $\displaystyle |\nabla K(x,y)| \leq A |x-y|^{-n-1}$ for some $A\geq0$, 
 \item the operator $\displaystyle f \mapsto \int_{\R^n} K(x,y) f(y)\, dy$ 
  initially defined on
 $C_0^\infty(\R^n)$ extends to a bounded operator on $L^2(\R^n)$.
\end{enumerate}
If $P:\R^n\to\R^n$ is a polynomial, then the operator $T$ defined by
\[
T[f](x) = \int_{\R^n} e^{i P(x,y)} K(x,y) f(y)\, dy
\]
can be extended to a bounded operator from $L^q(\R^n)$ to itself, with $1<q<\infty$. The
bound of this operator may depend on $K$, $q$, $n$ and
the degree $d$ of $P$ but is otherwise independent of the coefficients of $P$.
\end{Theorem}

\begin{proof}[Proof of Theorem \ref{thm:Lp bound}] 
The first step of the proof is to use 
the T(1) Theorem (p. 294 in \cite{St93}) on $S_\tau$. The T(1) Theorem 
says that if $S$ is a continuous linear
mapping from $\mathcal{S}$ to $\mathcal{S'}$ satisfying (\ref{it:xi-w}) and (\ref{it:z-z'}) (when $k=0$) and $S$ and
$S^*$ are restrictly bounded in the sense of (\ref{eqn:restr bound}), then $S$ extends to a bounded linear operator
from $L^2$ to itself. In our case, this means $S_\tau$ extends to a bounded linear operator.
However, since all of the constants in Lemma \ref{lem:Calderon-Zygmund operators} and 
Lemma \ref{lem:rest bound} are independent of $\tau$, it follows that $S_\tau$ is a bounded linear operator
from $L^2$ to itself with constants independent in $\tau$.

Next, $S_\tau$  satisfies the hypotheses of Theorem \ref{thm:ricci-stein}, so $\Ts$ is a bounded linear
operator from $L^q$ to itself for $1<q<\infty$ with a constant independent of $\tau$
but possibly depending on the $L^q$ constant of $S_\tau$ and the degree of $\tau T$ (which is
$\leq \deg p$), both of which are independent of $\tau$. 
\end{proof}

%
%
\section{Equivalence with NIS operators.}\label{sec:NIS equiv}

We now 
generate an OPF operator $\Ts$ 
from an NIS operator $\tT$ on a polynomial model 
$M^p$. Let $\tK(p,q)$ be the kernel of an NIS operator $\tT$.
On $\C\times\R$,
each kernel $\tK$ can be associated with a kernel $k$  by setting
\[
k(z,w,t-s) = \tK((z,t),(w,s)).
\]
The convolution structure in $t$ follows from the property that a polynomial model
is translation invariant in $t=\Rre z_2$.
Thus we have (for appropriate $\vp$),
\[
\tT[\vp](z,t) = \int_{\C\times\R} \tK((z,t),(w,s))\vp(w,s) \, dw ds 
= \int_{\C\times\R} k(z,w,t-s)\vp(w,s) \, dw ds.
\]
We set 
\begin{equation}\label{eqn:K def}
\Ks(z,w) = \int_{\R} e^{- i \tau t} k(z,w,t)\,dt
\end{equation}
and observe we also have
\[
k(z,w,t) = \frac {1}{2\pi}\int_\R e^{ i t\tau} \Ks(z,w)\,dt.
\]
The integrals representing $\Ks(z,w)$ and $k(z,w,t)$ do not necessarily converge. For a 
tempered distribution $T$
and a Schwartz function $\vp$, 
we know that if $\mathcal{F}$ represents the partial Fourier transform
in $t$, by definition, $\langle \mathcal{F} T, \vp\rangle = \langle T, \mathcal{F}\vp\rangle$.
As an integral, this corresponds to:
\begin{equation}\label{eqn:string of equalities}
\langle \mathcal{F} T, \vp\rangle
=\int_{\C\times\R}\hspace{-7.07pt} k(z,w,t)\int_{\R} e^{-i t \tau}\vp(w,\tau)\, d\tau dw dt
= \int_{\C\times\R} \int_{\R} k(z,w,t)e^{-i t \tau}\, dt\, \vp(w,\tau)\, dwdt.
\end{equation}
We make sense of \eqref{eqn:K def} by the string of equalities in
\eqref{eqn:string of equalities}, and we say
the integral $ \int_{\R} k(z,w,t)e^{-i t \tau}\, dt$ is defined
in the sense of Schwartz distributions. We similarly justify writing
$k(z,w,t) = \frac {1}{2\pi}\int_\R e^{ i t\tau} \Ks(z,w)\, d\tau$. If one of (or both of) the kernels is actually
in $L^1(\R)$ (in $t$ or $\tau$), then the integral defined in the sense of Schwartz distributions
agrees with the standard definition.

\subsection{An NIS Operator on $\C\times\R$ generates an OPF operator $\Ts$ on $\C$.}\label{subsec:up to down}

\n \begin{Theorem} \label{thm:up to down}
An NIS operator $\tT$ of order $m \leq 2$ on a polynomial model 
$M^p = \{(z_1,z_2)\in\C^2 : \Imm z_2 = p(z_1)\}$
generates an OPF operator $\Ts$ of order $m$ with respect to the polynomial $p$.
\end{Theorem}

\begin{Remark} \label{rem:pleasework}
The approximation conditions, (b) in the definition of OPF operators  and
(a) in the definition of NIS operators, 
imply one another since a partial Fourier transform is a 
continuous operator on the space of Schwartz distributions. Also, the adjoint conditions 
(f) from OPF operators and (e) from NIS operators,  allow us
to focus only $k$ and $\Ks$ as the computations will automatically apply to $k^*$ and $\Ks^*$.
\end{Remark}

Theorem \ref{thm:up to down} is proved in a series of lemmas. 
We first show that if $\tK$ is an NIS operator of order $m \leq 2$, then $\Ks$ is the kernel for a family
$\Ts$ of operators on $\C$.

The proof that $\Kse$ satisfies the size conditions \eqref{it:size} and \eqref{it:size2} is broken
into two lemmas. We handle the $m\leq 1$ case and the $m=2$ case.
\begin{Lemma}\label{lem:size up to down} 
If $m\leq 1$, the kernel $\Kse$ satisfies the size condition \eqref{it:size}.
\end{Lemma}

\begin{proof} It is enough to assume 
\[
\Ystp^J = \Ms^n = \epl \frac{\p^{n}}{\p\tau^{n}} \emi
\]
where $|J|=n$. 
Let $\eta\in C^\infty_c(\R)$ so that $\eta \equiv 1$ on $[-1,1]$, $0\leq\eta\leq 1$, and
$|\eta^{(n)}| \leq c_n$. Also, let $\eta_A(t) = \eta(t/A)$. We will estimate
\[
\frac{\p^n}{\p\tau^n} \int_{\R} e^{- i \tau(t + T(w,z))} \Ke(z,w,t) \eta_A(t)\, dt,
\]
and \eqref{it:size} will follow by sending $A\to\infty$.
The integral is compactly supported and the integrand is smooth, 
so we can apply the derivatives
inside of the integral. Integrating by parts $(n+k)$ times shows 
\begin{align*}
&c_n\left|\int_{\R} e^{- i \tau(t + T(w,z))} \big(t+T(w,z)\big)^n\Ke(z,w,t) \eta_A\big(t+T(w,z)\big)\, dt\right|\\
&= \frac{c_{n+k}}{|\tau|^{n+k}}\left|\int_{\R} e^{- i \tau(t + T(w,z))}
\frac{\p^{n+k}}{\p t^{n+k}}\Big(\big(t+T(w,z))^n \Ke(z,w,t) \eta_A(t+T(w,z))\Big)\, dt\right|\\
&= \frac{c_{n+k}}{|\tau|^{n+k}}\left|\int_{\R} \hspace{-2.51785pt}e^{- i \tau(t + T(w,z))}
\sum_{j=0}^{n+k} c_j \frac{\p^{j}}{\p t^{j}}\Big(\big(t+T(w,z))^n \Ke(z,w,t)\Big)
\eta_A^{(n+k-j)}\big(t+T(w,z)\big)\, dt\right|\\
&\leq \frac{c_{n+k}}{|\tau|^{n+k}} \sum_{j=1}^{n+k}\bigg[ \int_{|t+T(w,z)|\leq \Lambda(z,|w-z|)}
\hspace*{-1in}\Lambda(z,|w-z|)^{n-1-j}|w-z|^{m-2} \frac{1}{A^{n+k-j}}\, dt\\
&+ \int_{\Lambda(z,|w-z|) \leq |t+T(w,z)| \leq 2A}\hspace*{-1in}
|t+T(w,z)|^{n-1-j} \mu(z,|t+T(w,z)|)^{m-2}\frac{1}{A^{n+k-j}} \left|\eta^{(n+k-j)}\big(\tfrac{t+T(w,z)}A\big)
\right|\,dt \bigg].
\end{align*} 
If $j=n+k$, then 
\begin{align*}
& \frac{1}{|\tau|^{n+k}}\int_{|t+T(w,z)|\leq \Lambda(z,|w-z|)}
\hspace*{-1in}\Lambda(z,|w-z|)^{n-1-(n+k)}|w-z|^{m-2} \, dt\\
&+ \frac{1}{|\tau|^{n+k}}\int_{\Lambda(z,|w-z|) \leq |t+T(w,z)| \leq 2A} \hspace*{-1in}
|t+T(w,z)|^{n-1-j} \mu(z,|t+T(w,z)|)^{m-2}\frac{1}{A^{n+k-j}} \eta\big(\tfrac{t+T(w,z)}A\big)\,dt\\
&\leq c_{n+k}\frac{|w-z|^{m-2}}{|\tau|^{n+k}\Lambda(z,|w-z|)^k} 
+ \frac{|w-z|^{m-1}}{|\tau|^{n+k}} \int_{\Lambda(z,|w-z|) \leq |t+T(w,z)|} \hspace*{-1in}
|t+T(w,z)|^{-1-k} \mu(z,|t+T(w,z)|)^{-1}\, dt.
\end{align*} 
Using the substitution $s = \mu(z,|t+T(w,z)|)^{-1}$,  $|\frac{ds}{dt}| \sim
\frac{1}{\mu(z,|t+T(w,z)|)|t+T(w,z)|}$, so
\begin{align*}
\frac{|w-z|^{m-1}}{|\tau|^{n+k}}\int_{\Lambda(z,|w-z|) \leq |t+T(w,z)|}\hspace*{-1in}&
|t+T(w,z)|^{-1-k} \mu(z,|t+T(w,z)|)^{-1}\, dt \\
& \sim \frac{|w-z|^{m-1}}{|\tau|^{n+k}} \int_{|s|\leq \frac{1}{|w-z|}} \frac{1}{\Lambda(z,\frac1s)^k}\,ds
\leq c_{n+k}\frac{|w-z|^{m-2}}{|\tau|^{n+k}\Lambda(z,|w-z|)^k} .
\end{align*}
If $j< n+k$, then using the support condition of
$\eta_A^{(j)}\big(t+T(w,z)\big)$ that $|t+T(w,z)|\sim A$, the estimate simplifies to 
\begin{align*}
& \frac{1}{|\tau|^{n+k}}\int_{|t+T(w,z)|\leq \Lambda(z,|w-z|)}
\Lambda(z,|w-z|)^{n-1-j}|w-z|^{m-2} \frac{1}{A^{n+k-j}}\, dt\\
&+ \frac{1}{|\tau|^{n+k}}\int_{\Lambda(z,|w-z|) \leq |t+T(w,z)| \leq 2A} \hspace*{-1in}
|t+T(w,z)|^{n-1-j} \mu(z,|t+T(w,z)|)^{m-2}\frac{1}{A^{n+k-j}} \eta^{(n+k-j)}\big(\tfrac{t+T(w,z)}A\big)\,dt\\
&\leq  c_{n+k} \Lambda(z,|w-z|)^{n-j}|w-z|^{m-2} \frac{1}{A^{n+k-j}}
+ c_{n+k} A^{n-1-j} \mu(z,A)^{m-2}\frac{1}{A^{n+k-j+1}}
\stackrel{A\to\infty}{\longrightarrow} 0.
\end{align*}
This complete the proof for $m\leq 1$.
\end{proof}

\begin{Lemma}\label{lem:size up to down2} 
If $m=2$, the kernel $\Kse$ satisfies the size conditions \eqref{it:size} and \eqref{it:size2}.
\end{Lemma}

\begin{proof} As in Lemma \ref{lem:size up to down2}, we can assume that
\[
\Ystp^J = \Ms^n = \epl \frac{\p^{n}}{\p\tau^{n}} \emi
\]
where $|J|=n$. 

We first show the case $\mu(z, \frac 1\tau) \geq |w-z|$ and assume
$n=0$. From the definition of NIS operators,
$|\Ke(z,w,t)| \leq \frac{c_1}{\Lambda(z,|w-z|) + |t+T(w,z)|}$ and
$|\frac{ \p\Ke}{\p t}(z,w,t)| \leq \frac{c_2}{\Lambda(z,|w-z|)^2 + |t+T(w,z)|^2}$. 
Since $\Ke$ is not integrable on $\R$, we need to integrate by parts  to obtain 
an estimate on $\Kse$. However, since $|w-z|$ is small, we need to be careful to
integrate by parts as few times as possible and then only for large $t$. 
Let $A$ be a large number. 
\begin{align}
\Bigg| &\int_{|t+T(w,z)| \leq \frac{A}{|\tau|}} e^{- i \tau t}\Ke(z,w,t)\, dt \Bigg| \leq
 \left| \int_{|t+T(w,z)| \leq \Lambda(z,|w-z|)} e^{- i \tau t}\Ke(z,w,t)\, dt \right|\nn \\
&+ \left| \int_{\Lambda(z,|w-z|) \leq |t+T(w,z)| \leq \frac{1}{|\tau|}} e^{- i \tau t}\Ke(z,w,t)\, dt \right|
+ \left| \int_{\frac{1}{|\tau|}\leq |t+T(w,z)| \leq \frac{A}{|\tau|}} e^{- i \tau t}\Ke(z,w,t)\, dt \right|\nn\\
&\les 1 +  \int_{\Lambda(z,|w-z|) \leq |t+T(w,z)| \leq\frac{1}{|\tau|}} \frac{1}{|t+T(w,z)|}\, dt \nn\\
&+ \frac{1}{|\tau|} \left| \int_{\frac{1}{|\tau|}\leq |t+T(w,z)| \leq \frac{A}{|\tau|}}
e^{- i \tau t}\frac{\p\Ke}{\p t}(z,w,t) \, dt \right|
+ \frac{1}{|\tau||t+T(w,z)|}\bigg|_{|t+T(z,w)|=\frac 1\tau}^{|t+T(z,w)|=\frac A\tau} \label{eq:log calc}\\
&\les 1 + \log\left(\frac {1/|\tau|}{\Lambda(z,|w-z|)}\right) 
+ \frac{1}{|\tau|}\int_{\frac{1}{|\tau|}\leq |t+T(w,z)| \leq \frac{A}{|\tau|}} \frac{1}{\big(t+T(w,z)\big)^2}\, dt\nn \\
&\les 1 + \log\left(\frac {1/|\tau|}{\Lambda(z,|w-z|)}\right). \nn
\end{align}
This is actually the estimate we are looking for since
$\log \Big(\frac{1/|\tau|}{\Lambda(z,|w-z|)}\Big) \sim  \log\Big(\frac{\mu(z,\frac 1\tau)}{|w-z|}\Big)$.
Also, the estimate is independent of $A$, so we can let $A\to\infty$. 

Now
assume $k\geq 1$. Let $\eta\in C^\infty_c(\R)$, $0 \leq \eta \leq 1$, 
$\supp\eta(\cdot + T(w,z))\subset [-2,2]$,
$\eta\big(t+T(w,z)\big) =1$ if $|t|\leq 1$, and $\eta^{(k)}\big(t+T(w,z)\big) \leq c_k$. 
We show the case
$|w-z| \geq \mu(z,\frac 1\tau)$. Let $A\in \R$ be large. Integration by parts
$n+k$ times shows:
{\small
\begin{align*}
&\bigg |\frac{\p^n}{\p\tau^n} \int_\R e^{- i \tau(t+T(w,z))} 
\Ke(z,w,t) \eta(\tfrac{t+T(w,z)}A)\, dt\bigg |\\
&= \left|\sum_{j=0}^{n+k} \frac{c_j}{\tau^{n+k}} \int_\R e^{- i \tau(t+T(w,z))} 
\frac{\p^j}{\p t^j}\big(\big(t+T(w,z)\big)^n \Ke(z,w,t)\big) 
\frac{1}{A^{n+k-j}}\frac{d^{n+k-j}\eta}{dt^{n+k-j}}\Big(\tfrac{t+T(w,z)}A\Big) \, dt \right|\\
&\leq \frac C{|\tau|^{n+k}}\left( \sum_{j=0}^{n+k-1} A A^{n-1-j} A^{-n-k+j}
+ \int_\R \left| \frac{\p^{n+k} }{\p t^{n+k} }\Big(\big(t+T(w,z)\big)^n \Ke(z,w,t)\Big) \right |\, dt \right)\\
&\leq\frac C{|\tau|^{n+k}} \Bigg( \frac{1}{A^k} + \int_{|t+T(w,z)|\leq \Lambda(z,|w-z|)}
\hspace{-.5in}\Lambda(z,|w-z|)^{-(k+1)}\, dt
+ \int_{|t+T(w,z)|\geq \Lambda(z,|w-z|)} \hspace{-.5in}|t+T(w,z)|^{-(k+1)}\, dt \Bigg)\\
&\leq \frac C{|\tau|^{n+k}}\left( \frac{1}{A^{k}} + \frac{1}{\Lambda(z,|w-z|)^k}\right).
\end{align*} }
Sending $A\to\infty$ yields the desired estimate.

We have one estimate left to compute: the case $|w-z|<\mu(z,\frac 1\tau)$ and $n\geq 1$.
Let $A$ be a large number. Let $0 \leq \psi_1,\psi_2^A \leq 1$ so that $1 = \psi_1+\psi_2^A$ on
$[-A,A]$. Let
$\supp \psi_1 \subset [-2,2]$ and $\supp \psi_2^A \subset \{ t : |t| \in [\frac 32, 2A]\}$, and
assume $|\frac{\p^n}{\p t^n}\psi_2^A| \leq \frac{c_n}{A^n}$ if $|t| \geq \frac{A}{2}$
and $|\frac{\p^n\psi_1}{\p t^n}|, |\frac{\p^n \psi_2^A}{\p t^n}| \leq c_n$ if $|t| \leq 2$.
Since $|z-w| \leq \mu(z, \frac 1\tau)$, $\Lambda(z,|z-w|) \les \frac 1\tau$. 
\begin{align*}
&\bigg |\frac{\p^n}{\p\tau^n} \int_\R e^{- i \tau(t+T(w,z))} 
\Ke(z,w,t) \Big(\psi_1\big(\tau(t+T(w,z))\big) + \psi_2\big(\tau(t+T(w,z))\big)\Big)\, dt\bigg |\\
&\leq c_n\int_{|t+T(w,z)| \leq \frac 2{|\tau|}} |t+T(w,z)|^n |\Ke(z,w,t)|\, dt \\
&+ \sum_{j=0}^n c_j \left| \int_{\R} \big(t+T(w,z)\big)^n \frac{\p^j \psi_2^A\big(\tau(t+T(w,z))\big)}{\p \tau^j}
\Ke(z,w,t) e^{- i \tau(t+T(w,z))}\, dt \right|
\end{align*}
Picking an arbitrary term and integrating by parts $(n+2)$ times, we have {\small
\begin{align*}
&\left| \int_{\R} \big(t+T(w,z)\big)^n \frac{\p^j \psi_2^A(\tau(t+T(w,z)))}{\p \tau^j}
\Ke(z,w,t) e^{- i \tau(t+T(w,z))}\, dt \right| \\
&\leq c_{n+2} \sum_{k=0}^{n+2}\int_{\R}\left| \frac{1}{(t+T(w,z))^{n+2}}
\frac{\p^k}{\p t^k}\Big((t+T(w,z))^n \Ke(z,w,t)\Big)
\frac{\p^{n+2+j-k} \psi_2^A(\tau(t+T(w,z)))} { \tau^{-n-2+k}\,\p \tau^{j}\p t^{n+2-k}} \right| \, dt 
\end{align*} }
If $n+2+j-k \geq 1$, the term in the sum has support near $\frac 1{|\tau|}$ and
$\frac{A}{|\tau|}$, so it is bounded by
\begin{align*}
&\int_{\R}\left| \frac{\tau^{n+2-k}}{\tau^{n+2}}
\frac{\p^k}{\p t^k}\Big((t+T(w,z))^n \Ke(z,w,t)\Big) 
\frac{\p^{n+2+j-k} \psi_2^A(\tau(t+T(w,z)))} {\p \tau^{j}\p t^{n+2-k}} \right| \, dt \\
&\leq\frac{c_n}{|\tau|^{n+2}}\frac{1}{|\tau|^{n-1-k}} |\tau|^{n+2-k} \frac{1}{|\tau|}
+\frac{c_n}{|\tau|^{n+2}} \frac{A^{n-1-k}}{|\tau|^{n-1-k}} \frac{|\tau|^{n+2-k}}{A^{n+2-k+j}}
\frac{A}{|\tau|}  \stackrel{A\to\infty}{\longrightarrow} \frac{c_n}{|\tau|^n}.
\end{align*}
Finally, if $n+2+j-k=0$, then $j=0$ and $k=n+2$ and we have the estimate
\begin{align*}
\int_{\R}&\left| \frac{1}{\tau^{n+2}}
\frac{\p^{n+2}}{\p t^{n+2}}\Big((t+T(w,z))^n \Ke(z,w,t)\Big)
\psi_2^A\big(\tau(t+T(w,z))\big)\right| \, dt \\
&\leq \frac{c_n}{|\tau|^{n+2}} \int_{|t+T(w,z)|\geq \frac{1}{2|\tau|}}
\frac{1}{|t+T(w,z)|^3}\, dt = \frac{c_n}{|\tau|^n}.
\end{align*}
\end{proof}

\begin{Lemma}\label{lem:cancel up to down} 
The operator $\Ts$ has the $w$-cancellation condition \eqref{it:cancel w}.
\end{Lemma}

\begin{proof} Let $\Ystp^{J}$ be a product of $|J|$ operators of the form
$\Ystp^j = \Zstp,\ \Zbstp,\ \Ms$ where $|J| = \ell +n$ and
$n = \#\{j : \Ystp^j = \Ms\}$ and let $\vp \in C^{\infty}(D(z_0,\delta))$. We have
\[
\Kse(z,w) = \int_\R e^{- i \tau t} \Ke(z,w,t)\,dt ,
\]
so that integration by parts yields
\begin{align*}
\Zstp \Kse(z,w) &=\Zstp\int_\R e^{- i \tau t}\Ke(z,w,t) \, dt \\
&= \frac {\p}{\p z} \int_\R e^{- i \tau t}\Ke(z,w,t) \, dt - \int_\R  \tau \frac{\p p}{\p z}(z)
e^{- i \tau t}\Ke(z,w,t) \, dt \\
&= \int_\R  e^{- i \tau t} L\Ke(z,w,t)\, dt.
\end{align*}
Similarly, $\Zbstpz \Kse(z,w)= \int_\R  e^{- i \tau t} \LL_z\Ke(z,w,t)\, dt$.
Also, recalling that $\MM f(z,w) = - i \big(t+T(w,z)\big)f(z,w)$, we have
$\Ms\Kse(z,w) =   \int_\R  e^{- i \tau (t+T(w,z))} \MM\Ke(z,w,t)\, dt$.
Thus,
\[
\int_\C \Ystp^J\Kse(z,w)\vp(w)\, dw
= \int_\C \int_\R e^{- i \tau t} \YY^Jk(z,w,t) \vp(w)\, dt dw,
\]
with the correspondence that if $\Ystp^j = \Zstp,\Zbstp, \Ms$, then
$\YY^j = L, \LL, \MM$ respectively.
Integrating $(n+k)$ times gives us:
\begin{align}
\int_{\C} \Ystp^J &\Kse(z,w)\vp(w)\, dw
= \iint_{\C\times\R} (\YY^J \Ke)(z,w, t) e^{- i \tau t}\vp(w)\, dt dw \nn\\
&=  \frac{c_{n+k}}{\tau^{n+k}} \iint_{\C\times\R} 
\left(\frac{\p^{n+k}}{\p t^{n+k}}\YY^J\right)\Ke(z,w, t)
e^{- i \tau t}\vp(w) \eta(w,t) \, dt dw \nn \\
&\ \ + \frac{c_{n+k}}{\tau^{n+k}} \iint_{\C\times\R} 
\left(\frac{\p^{n+k}}{\p t^{n+k}}\YY^J\right)\Ke(z,w, t)
e^{- i \tau t}\vp(w)(1-\eta(w,t))\, dt dw \label{eq:w cancel int by parts}
\end{align}
where $\eta\in C^\infty_c(\C\times\R)$ is a bump function on 
$B_{NI}((z,0),\delta)$.
To estimate the integrals in \eqref{eq:w cancel int by parts}, the strategy is to 
expand $\left(\frac{\p^{n+k}}{\p t^{n+k}}\YY^J\right)\Ke(z,w, t)$
and estimate an arbitrary term.
It is important to remember that in $\YY^J$, $n$ of the terms are $\MM$ and an
$L$ or $\LL$ can hit either an $\MM$ term or $\Ke(z,w,t)$.

Expanding $\left(\frac{\p^{n+k}}{\p t^{n+k}}\YY^J\right)\Ke(z,w, t)$, we see
\begin{align}
&\frac{\p^{n+k}}{\p t^{n+k}}\YY^J \Ke(z,w, t) \nn \\
&= \frac{\p^{n+k}}{\p t^{n+k}} \left[\sum_{|J_0| + \cdots + |J_n|= \ell} \left(
c_{|J_0|, \ldots, |J_n|}
\XX^{J_0} \Ke(z,w, t) \prod_{j=1}^n(- i ) \XX^{J_j}\big(t+T(w,z)\big)\right)\right]\nn \\
&=  \sum_{\atopp{|J_0| + \cdots + |J_n|= \ell} {\ell_0 + \cdots + \ell_n = {n+k}}}
c_{|J_0|, \ldots, |J_n|}  c_{\ell_0,\ldots, \ell_n} 
\frac{\p^{\ell_0}}{\p t^{\ell_0}} \XX^{J_0} \Ke(z,w, t) \prod_{j=1}^n 
\frac{\p^{\ell_j}}{\p t^{\ell_j}} \XX^{J_j}\big(t+T(w,z)\big), \label{eq:Kep derivs}
\end{align}
\n where $\XX^{J_j}$ is an operator composed only of $\XX^j = L$ and $\LL$. 
We pick an arbitrary term from the sum and show that it has the desired bound.
Taking an arbitrary term from \eqref{eq:Kep derivs}, we estimate the integrals
from \eqref{eq:w cancel int by parts} which reduce to the following two integrals:
\[
I = \left| \frac{1}{\tau^{n+k}} \iint_{\C\times\R} 
\frac{\p^{\ell_0}}{\p t^{\ell_0}} \XX^{J_0} \Ke(z,w, t) \prod_{j=1}^n 
\frac{\p^{\ell_j}}{\p t^{\ell_j}} \XX^{J_j}\big(t+T(w,z)\big)
e^{- i \tau t}\vp(w) \eta(w,t) \, dt dw \right|
\]
and
{\small \[ 
II = \left | \frac{1}{\tau^{n+k}} \iint_{\C\times\R} 
\frac{\p^{\ell_0}}{\p t^{\ell_0}} \XX^{J_0} \Ke(z,w,t)  
\prod_{j=1}^n \frac{\p^{\ell_j}}{\p t^{\ell_j}} \XX^{J_j}\big(t+T(w,z)\big)
e^{- i \tau t}\vp(w)(1-\eta(w,t))\, dt dw \right|
\]}
where $|J_0| + \cdots + |J_\ell|= \ell$ and $\ell_0 + \cdots + \ell_n = {n+k}$.
Using Proposition  \ref{prop:t+T deriv} and the 
cancellation condition \eqref{it:cancel w}, $I$ has the estimate:

\begin{align*}
& I \leq \frac{c_{|J_0|,\ell_0}}{|\tau|^{n+k}} \frac{\delta^{m-|J_0|}}{\Lambda(z,\delta)^{\ell_0}} 
\hspace{-1pt}
\sup_{(w,t)} \hspace{-4pt}\sum_{|I|\leq N_{|J_0|,\ell_0}} \hspace{-15pt} \delta^{|I|} \left| \XX^I\Big(
e^{- i \tau t}\vp(w) \prod_{j=1}^n \left(\frac{\p^{\ell_j}}{\p t^{\ell_j}} \XX^{J_j}\big(t+T(w,z)\big)\right)
\eta(w,t)\Big)\right| \\
&\leq \frac{c_{|J_0|,\ell_0}}{|\tau|^{n+k}} \frac{\delta^{m-|J_0|}}{\Lambda(z,\delta)^{\ell_0}} 
\sup_{(w,t)} \sum_{|I|\leq N_{|J_0|,\ell_0}} \delta^{|I|} 
\sum_{|I_0| + \cdots + |I_{n+1}|=|I|} c_{I_0,\ldots, I_{n+1}} \Bigg| \XX^{I_0}\big(
e^{- i \tau t}\vp(w)\big) \\
&\hspace{171pt} \times\prod_{j=1}^n \left( \XX^{I_j}\frac{\p^{\ell_j}}{\p t^{\ell_j}} 
\XX^{J_j}\big(t+T(w,z)\big)\right) \XX^{I_{n+1}}\eta(w,t) \Bigg| \\
&\leq  \frac{c_{n,\ell,k}}{|\tau|^{n+k}}\Lambda(z,\delta)^{-k} \delta^{m-\ell}
\sup_{(w,t)}\sum_{|I_0|\leq N_{|J_0|,\ell}} \delta^{|I_0|} \big|\XX^{I_0}(e^{- i \tau t}\vp(w))\big| \\
&=  \frac{c_{n,\ell,k}}{|\tau|^{n+k}}\Lambda(z,\delta)^{-k} \delta^{m-\ell}
\sup_{w}\sum_{|I_0|\leq N_{|J_0|,\ell}} \delta^{|I_0|} \big|X^{I_0}_\tau\vp(w)\big|. 
\end{align*} 
To estimate $II$, we use size estimates and the support size of $\vp$. 
{\small
\begin{align}
II &\leq \frac{c_{n,\ell} \|\vp\|_{L^\infty}}{|\tau|^{n+k}} \int_{|w-z_0|\leq\delta}
\int_{|t+T(w,z)|\geq \Lambda(z,\delta)}
\frac{\dni(z,w,t)^{m-2-|J_0|}}{\Lambda(z,\dni(z,w,t))^{1+\ell_0}} \nn\\
&\hspace{152pt}
\times\frac{\Lambda(z,\dni(z,w,t))^n}{\dni(z,w,t)^{|J_1|+\cdots |J_n|} 
\Lambda(z,\dni(z,w,t))^{\ell_1 +\cdots+ \ell_n}}  \, dt dw \nn\\
&\leq  \frac{c_{n,\ell}}{|\tau|^{n+k}} \|\vp\|_{L^\infty} \int_{|w-z_0|\leq\delta}
\int_{|t+T(w,z)|\geq \Lambda(z,\delta)}\hspace{-45pt} \mu(z,t+T(w,z))^{m-\ell-2} 
\frac{1}{|t+T(w,z)|^{n+k-n+1}} \, dt dw. \label{eq:w cancel m=2}
\end{align}
If $m\leq2$ or $m=2$ and $\ell\geq1$, then we use the substitution
$s= \mu\big(z,t+T(w,z)\big)^{-1}$, so $|\frac 1s \frac{ds}{dt}| \sim |t+T(w,z)|^{-1}$ and
\eqref{eq:w cancel m=2} becomes
\begin{align*}
II &\leq \frac{c_{n,\ell}}{|\tau|^{n+k}}\|\vp\|_{L^\infty}\Lambda(z,\delta)^{-k} \delta^2
\int_{|s|\leq \frac 1\delta} s^{1-m+\ell} \, ds 
\leq \frac{c_{n,\ell}}{|\tau|^{n+k}}\|\vp\|_{L^\infty}\Lambda(z,\delta)^{-k} 
\delta^{m-\ell}.
\end{align*} }
If $m=2$, $\ell=0$, and $k\geq 1$, then a straightforward integration shows that
$II \leq \frac{c_{n,\ell}}{|\tau|^{n+k}}\|\vp\|_{L^\infty}\Lambda(z,\delta)^{-k} 
\delta^{2}$.
The integral in \eqref{eq:w cancel m=2} 
diverges if $m=2$ and $\ell=k=0$, so we must estimate the tail term more carefully in
this case. With $m=2$, $\ell=0$, and $k=0$, \eqref{eq:w cancel m=2} simplifies to
\[
II \leq \frac{1}{|\tau|^n} \left|\int_{\C}\int_{\R} e^{-i t\tau} \frac{\p^{\ell_0}}{\p t^{\ell_0}} \Ke(z,w,t)
\big(t+T(w,z)\big)^{n-(n-\ell_0)} \vp(w)\big(1-\eta(w,t)\big)\, dt dw \right|. 
\]
The key to this estimate is to recognize that 
$\frac{\p^{\ell_0}}{\p t^{\ell_0}} \Ke(z,w,t)\big(t+T(w,z)\big)^{\ell_0}$ satisfies the estimates of
an order 2 NIS operator. To integrate in $t$, we use
the argument of  (\ref{eq:log calc}) with $\delta$ replacing $|z-w|$
and see that 
\[
|II| \leq \frac{c_{n,0}}{|\tau|^n}     
\int_\C |\vp(w)| \log(\tfrac 1{\tau \Lambda(z,\delta)}) \, dw
\les \frac{c_{n,0}}{|\tau|^n}  \delta^2 \log(\tfrac{\mu(z,\frac 1\tau)}{\delta}) \|\vp\|_{L^\infty(\C)}.
\]
Note that $\log(\tfrac 1{\tau \Lambda(z,\delta)}) \sim \log(\tfrac{\mu(z,\frac 1\tau)}{\delta})$.
While this estimate is true for all $\tau$ and $\delta$, the previous
estimate of $II$ shows that we only have to consider the case
when $\delta \leq \mu(z,\frac 1\tau)$ or equivalently, $\tau\Lambda(z,\delta) \leq 1$.
\end{proof}

\begin{Lemma}\label{lem:up to down cancel tau} 
The kernel $\Kse$ satisfies the $\tau$-cancellation condition \eqref{it:cancel tau}.
\end{Lemma}

\begin{proof} 
Since $\mathcal{F}^{-1} \mathcal{F} = I$ in the sense of Schwartz distributions, 
\[
|\XX^Jk(z,w,t)| \leq C_{|J|} \frac{\mu(z, t+T(w,z))^{m-|J|}}{V(z,\mu(z,t+T(w,z))))}
\]
implies $\frac1{2\pi}\int_\R \Xstp^J \Big(e^{ i \tau t} \Kse(z,w)\Big)\, d\tau 
= \XX^Jk(z,w,t)$ satisfies the same estimates.
\end{proof}

The proof of  Theorem \ref{thm:up to down} is complete.

\subsection{An OPF operator $\Ts$ on $\C$ generates an NIS operator $\tK$ on $\C\times\R$.}
\label{subsec:down to up}

\n\begin{Theorem} \label{thm:down to up}
A OPF operator $\Ts$ of order $m\leq 2$ with respect to the subharmonic, nonharmonic
polynomial $p$ generates
an NIS operator $\tK$ of order $m \leq 2$
on the polynomial model  $M^p = \{(z_1,z_2)\in\C^2: \Imm z_2 = p(z_1)\}$.
\end{Theorem}
 
We prove Theorem \ref{thm:down to up} in the same manner as 
Theorem \ref{thm:up to down}. Remark \ref{rem:pleasework}
applies to Theorem \ref{thm:down to up} as well.

\begin{Lemma} \label{lem: NIS cancel} 
The operator $\tK$ satisfies the NIS cancellation conditions \eqref{eq:NIS cancel}.
\end{Lemma}

\begin{proof} Let $\vp\in C^\infty_c\big(B((z,t),\delta)\big)$. Also,
let $\vps(z,\tau) = \int_{\R} e^{-i \tau t} \vp(z,t)\,dt$ be the partial Fourier transform
in $t$ of $\vp(z,t)$.
Let $\eta\in C^\infty_c(\R)$ with
$\supp\eta \subset [-\frac 2{\Lambda(z,\delta)}, \frac 2{\Lambda(z,\delta)}]$ and $\eta(\tau) =1$
when $\tau\in [-\frac 1{\Lambda(z,\delta)}, \frac 1{\Lambda(z,\delta)}]$.
Let $\XX^J$ be a product of $|J|$ operators of the form of $\XX^j = \LL_z$ and $L_z$. Then
\begin{align*}
&\XX^J\iint_{\C\times\R} \hspace{-9pt}\Ke(z,w, t-s)\vp(w, s)\,dw ds
= \frac{1}{2\pi}\int_\C\int_\R\int_\R \XX^J \big(e^{ i \tau(t-s)} \Kse(z,w)\big)\vp(w,s)\, d\tau ds d w \\
&\hspace{-3pt}=\hspace{-3pt}
\frac{1}{2\pi} \int_{\C}\int_\R \hspace{-3pt}e^{ i t\tau} \Xstp^J \Kse(z,w)\vps(w,\tau)\, dw d\tau 
\hspace{-1.5pt}=\hspace{-1.5pt} 
\frac{1}{2\pi}\int_\R e^{ i t\tau} \Xstp^J \int_\C \hspace{-3pt}\Kse(z,w) \vps(w,\tau)\, dw \,\eta(\tau) d\tau \\
&+\frac{1}{2\pi} \int_\R e^{ i t\tau} \Xstp^J \int_\C \Kse(z,w) \vps(w,\tau)\, dw\,(1- \eta(\tau)) d\tau
= I + II.
\end{align*}
We estimate $I$ and $II$ separately. We first do the case $m\leq 1$ or $m=2$ and $|J|\geq 1$.
By \eqref{it:cancel w},
\begin{align*}
|I| &\leq c_{|J|} \delta^{m-|J|} \int_{\R} \sup_{w} \sum_{|I|\leq N_{|J|}} \delta^{|I|}
\Big| \Xstp^I\Big( \vps(w,\tau)\eta(\tau)\Big)\Big| \, d\tau \\
&\leq c_{|J|} \delta^{m-|J|}  \int_{\R} |\eta(\tau)| \sup_{w} \sum_{|I|\leq N_{|J|}} \delta^{|I|}
\|\XX^I \vp\|_{L^{1}(t)} \, d\tau \\
&\leq c_{|J|} \delta^{m-|J|} \frac{1}{\Lambda(z,\delta)}\sum_{|I|\leq N_{|J|}} \delta^{|I|}
\|\XX^I \vp\|_{L^{\infty}(\C\times\R)} \Lambda(z,\delta).
\end{align*}
The last line follows from H\"older's inequality and the size of $\supp\vp$. 
The only difference between the $m=2$, $J=0$ case and the previous estimate is the
logarithm term in \eqref{it:cancel w}. The term to estimate is
\begin{equation}\label{eq:m=2 cancel upstairs log estimate}
\Big| \int_\R \eta(\tau) \delta^2 \log(\tfrac 1{\tau\Lambda(z,\delta)})  \sup_w |\vps(w,\tau)|\, d\tau  \Big|
\end{equation}
However, integration shows that $\int_0^{\frac{1}{\Lambda(z,\delta)}} \log(\frac{1}{\tau\Lambda(z,\delta)})\,
d\tau = \frac1{\Lambda(z,\delta)}$, so \eqref{eq:m=2 cancel upstairs log estimate} simplifies to
\[
\delta^2 \|\vps\|_{L^\infty(\C\times\R)}
 \int_0^{\frac1{\Lambda(z,\delta)}}  \log(\tfrac 1{\tau\Lambda(z,\delta)}) \, d\tau  
= \delta^2 \|\vps\|_{L^\infty(\C\times\R)} \frac{1}{\Lambda(z,\delta)} \les  \delta^2 \|\vp\|_{L^\infty(\C\times\R)}.
\]
We estimate $II$ in a similar fashion. We first cover the case when
$m\leq 1$ or $m=2$ and $|J|\geq 1$.
\begin{align}
|II| &=\frac{1}{2\pi}
 \left|\int_{\R}e^{ it\tau}\big(1-\eta(\tau)\big) \frac{1}{\tau^2}\left(\Xstp^J \int_{\C} \tau^2 \Kse(z,w)
\vps(w,\tau) \, dw\right) d\tau\right| \nn\\
&\leq c_{|J|} \int_{|\tau|>\frac{1}{\Lambda(z,\delta)}} |\tau|^{-2} \delta^{m-|J|}
\sum_{|I|\leq N_{|J|}} \delta^{|I|} \| \tau^2 \Xstp^I \vps(w,\tau)\|_{L^\infty(w)}\, d\tau.
\label{eq:II sum}
\end{align}
The terms in the sum can be rewritten the more useful way:
\begin{align}
\| \tau^2 X^I \vps(w,\tau)\|_{L^\infty(w)}
&= \sup_{w} \Big|\frac{1}{2\pi}\int_{\R} \tau^2 \Xstp^I e^{ i \tau t} \vp(w,t)\, dt \Big| \nn\\
&\hspace{-5pt}= c \sup_{w} \left|\int_{\R}  e^{ i \tau t}\left(\frac{\p^2}{\p t^2}\XX^I\vp(w,t)\right)\, dt \right| 
\leq c_2  \Lambda(z,t)\left\| \frac{\p^2}{\p t^2} \XX^I \vp \right\|_{L^{\infty}(\C\times\R)}.
\label{eq:II sum Linfty bound}
\end{align}
Using the estimate from \eqref{eq:II sum Linfty bound} in \eqref{eq:II sum}, 
\begin{align}
|II| &\leq c_{|J|}\delta^{m-|J|} \int_{|\tau|>\frac{1}{\Lambda(z,\delta)}} |\tau|^{-2}
\sum_{|I|\leq N_{|J|}} \delta^{|I|} \left \| \left (\frac{\p^2}{\p t^2} \XX^I\right) 
\vp(w,t)\right\|_{L^\infty(\C\times\R)}\Lambda(z,\delta)\, d\tau\nn \\
&\leq c_{|J|}\delta^{m-|J|}\sum_{|I|\leq N_{|J|}} \delta^{|I|} \Lambda(z,\delta)^2
\left \| \left (\frac{\p^2}{\p t^2} \XX^I\right)  \vp(w,t)\right\|_{L^\infty(\C\times\R)} 
\label{eqn:log doesn't screw things up}\\
&\leq c_{|J|}\delta^{m-|J|}\sum_{|I|\leq N_{|J|}'} \delta^{|I|}
\left \|\XX^I  \vp(w,t)\right\|_{L^\infty(\C\times\R)}.\nn
\end{align}
In the final estimate, we used the fact that $\Lambda(z,\delta)\frac{\p}{\p t}$ can be generated
by commutators of $\delta X$ terms.
As in $I$, the difference between the $m=2$, $J=0$ and the case already estimated
is the logarithm term in \eqref{it:cancel w}. However,
$\int_{\Lambda(z,\delta)^{-1}}^\infty \frac{|\log(\frac 1{\tau\Lambda(z,\delta)})|}{\tau^2}\, d\tau
= \Lambda(z,\delta)$, so we can repeat the estimate in \eqref{eqn:log doesn't screw things up}
replacing $|\tau|^{-2}$ with $\frac{|\log(\frac 1{\tau\Lambda(z,\delta)})|}{\tau^2}$ and achieve
the same conclusion.
\end{proof}

\begin{Lemma} \label{lem: NIS size} 
The operator $\tK$ has the NIS size conditions \eqref{eq:NIS size}.
\end{Lemma}

\begin{proof} It is enough to find the estimate on $|\Ke(z,w,t)|$.
We handle the $m=2$ separately. First assume $m\leq 1$. If
$\dni(z,w,t) = |z-w|$, then we break the integral in two pieces and estimate each piece separately.
\[
\int_{\R} e^{ i\tau t} \Kse(z,w)\,d\tau 
=\frac{1}{2\pi} \int_{|\tau|\leq \frac{1}{\Lambda(z,|w-z|)}}\hspace{-9pt}
e^{ i\tau t} \Kse(z,w)\,d\tau 
+ \frac{1}{2\pi}\int_{|\tau|\geq \frac{1}{\Lambda(z,|w-z|)}} \hspace{-9pt}
e^{ i\tau t} \Kse(z,w)\,d\tau.
\]
Estimating the first integral gives us:
\[
\left| \int_{|\tau|\leq \frac{1}{\Lambda(z,|w-z|)}} e^{ i\tau t} \Kse(z,w)\,d\tau  \right|
\leq c_0 \frac{|w-z|^m}{|w-z|^2\Lambda(z,|w-z|)} = c_0 \frac{\dni(z,w,t)^m}{V(z,\dni(z,w,t))}.
\]
The tail term is no harder: by \eqref{it:size} with $\ell=n=0$ and $k=2$,
\begin{align*}
\left| \int_{|\tau|\geq \frac{1}{\Lambda(z,|w-z|)}} e^{ i\tau t} \Kse(z,w)\,d\tau  \right|
&\leq c_2 \frac{|w-z|^m}{|w-z|^2\Lambda(z,|w-z|)^2} 
\int_{|\tau| \geq  \frac{1}{\Lambda(z,|w-z|)}} \frac{1}{\tau^2}\,d\tau\\
&\leq c_2 \frac{|w-z|^m}{|w-z|^2\Lambda(z,|w-z|)}.
\end{align*}
The case $\dni(z,w,t)=\mu(z,t+T(w,z))$ is the $\tau$-cancellation condition \eqref{it:cancel tau}.

Now assume $m=2$. The estimate to prove is
\[
|\Ke(z,w,t)| \leq C \frac{\dni(z,w,t)^2}{V(z,\dni(z,w,t))} = C \frac{1}{\Lambda(z,\dni(z,w,t))}.
\]
Let $\eta \in C^\infty_c(\R)$ where $\supp \eta \subset [-2,2]$,
$\eta(\tau) =1$ if $|\tau|\leq 1$, $0 \leq \eta \leq 1$,
and $\left|\frac{\p^k\eta}{\p \tau^k}(\tau)\right| \leq C_k$.
Let $\Lambda = \Lambda(z,\dni(z,w,t))$. We have
\[
\Ke(z,w,t) 
= \int_\R e^{ i\tau t}\Kse(z,w)\eta(\tau \Lambda)\,d\tau
+ \int_\R e^{ i\tau t}\Kse(z,w)(1-\eta(\tau \Lambda))\,d\tau  = I + II.
\] 
Before we estimate $I$, observe 
$\int_{\delta}^\infty \frac{\log s}{s^k}\, ds = -k \frac{\log s}{s^{k+1}}
+ \frac{k}{k+1} \frac{1}{s^k}$.
Also, with the change of variables $s = \frac{2\mu(z,\tfrac 1\tau)}{|w-z|}$,
$|\frac{\p s}{\p \tau}| \sim \frac{\mu(z,\frac 1\tau)}{|w-z|}\frac 1{|\tau|}$ and
$\Lambda(s,|w-z|) \sim \frac 1{|\tau|}$, so
\begin{align*}
I &\les \int_{|\tau| \leq \frac{2}{\Lambda}}
\log\left(\tfrac{2\mu(z,\tfrac 1\tau)}{|w-z|}\right) \, d\tau 
\sim \int_{|s| \geq \frac {\mu(z,\Lambda)}{|w-z|}} \frac{\log s}{s \Lambda(z,|w-z|s)}\, ds \\
&\sim \int_{\frac {\mu(z,\Lambda)}{|w-z|}}^\infty \inf_{j,k\geq 1} 
\frac{1}{|\A{jk}z| |w-z|^{j+k}} \frac{\log s}{s^{j+k+1}} \, ds\\
&\les \inf_{j,k\geq 1} \frac{1}{|\A{jk}z| |w-z|^{j+k}}
\left(\frac{\log\left(\frac {\mu(z,\Lambda)}{|w-z|}\right)}
{\left(\frac {\mu(z,\Lambda)}{|w-z|}\right)^{j+k+1}} 
+ \frac{1}{ \left(\frac {\mu(z,\Lambda)}{|w-z|}\right)^{j+k}} \right)\\
&\les \inf_{j,k\geq 1}\frac{1}{|\A{jk}z| |w-z|^{j+k}} 
\frac{|w-z|^{j+k}}{\mu(z,\Lambda)^{j+k}} 
\sim \frac{1}{\Lambda(z,\mu(z,\Lambda))} = \frac{1}{\Lambda}.
\end{align*}
To estimate $II$, we need to separate the cases $\Lambda = \Lambda(z,|w-z|)$ and
$\Lambda =  |t+T(w,z)|$. We first do the case $\Lambda = \Lambda(z,|w-z|)$.
By \eqref{it:size} with $k=2$ and $\ell=n=0$,
\[
II \les \int_{\frac{1}{\Lambda}}^\infty \frac{1}{\tau^2 \Lambda^2}\, d\tau
 \sim \frac{1}{\Lambda}.
\]
Now assume $\Lambda = |t+T(w,z)|$. Then {\small
\[
II \les \frac{1}{\big(t+T(w,z)\big)^2}\left| \int_{|\tau| \geq \frac{1}{|t+T(w,z)|}}\hspace{-0.67996pt}
\hspace*{-28pt}e^{ i \tau(t+T(w,z))} \frac{\p^2}{\p\tau^2}\left(
\emi\Kse(z,w)(1-\eta(\tau|t+T(w,z)|))\right) \,d\tau \right|.
\]  }
If both $\tau$-derivatives are applied to $\Kse$,
\begin{align*}
\frac{1}{\big(t+T(w,z)\big)^2}&\int_{|\tau| \geq \frac{1}{|t+T(w,z)|}}\left|
\frac{\p^2}{\p\tau^2}\Big(\emi\Kse(z,w)\Big)\right| \big(1-\eta(\tau|t+T(w,z)|)\big) \,d\tau \\
&\sim \frac{1}{\big(t+T(w,z)\big)^2} \int_{|\tau| \geq \frac{1}{|t+T(w,z)|}} \frac{1}{\tau^2} \, d\tau
\sim \frac{1}{\big(t+T(w,z)\big)}.
\end{align*}
Next, if one $\tau$-derivative is applied to $\Kse$ and one to $\eta$, then
\begin{align*}
\frac{1}{\big(t+T(w,z)\big)}&\int_{|\tau| \geq \frac{1}{|t+T(w,z)|}}\left|
\frac{\p}{\p\tau}\Big(\emi\Kse(z,w)\Big)\right| \eta'(\tau|t+T(w,z)|)) \,d\tau \\
&\sim \frac{1}{\big(t+T(w,z)\big)} \int_{|\tau| \sim \frac{1}{|t+T(w,z)|}} \frac{1}{\tau} \, d\tau
\sim \frac{1}{\big(t+T(w,z)\big)}.
\end{align*}
Finally, if $\eta$ receives both $\tau$-derivatives,
$\int_{|\tau| \geq \frac{1}{|t+T(w,z)|}}
\left| \Kse(z,w)\right| \eta''(\tau|t+T(w,z)|)) \,d\tau 
\sim \frac{1}{\big(t+T(w,z)\big)}$.

\end{proof}

Proving Theorem \ref{thm:up to down}
and Theorem \ref{thm:down to up} proves Theorem \ref{thm:correspondence}.


\section*{References.}


{\footnotesize
\centerline{\rule{9pc}{.01in}}
\bigskip
\bigskip
\centerline{Texas A\&M University}
\centerline{e-mail: andrew.raich@math.tamu.edu}
}

\end{document}